\documentclass{amsart}
\usepackage{amssymb}
\usepackage{graphics}
\usepackage{amscd}
\usepackage{enumerate}
\newtheorem{thm}{Theorem}[section]
\newtheorem{cor}[thm]{Corollary}
\newtheorem{lem}[thm]{Lemma}
\newtheorem{prop}[thm]{Proposition}

\newtheorem{claim}[thm]{Claim}
\theoremstyle{definition}
\newtheorem{defn}[thm]{Definition}

\newtheorem{qn}[thm]{Question}
\theoremstyle{remark}
\newtheorem{rem}[thm]{Remark}

\numberwithin{equation}{section}
\newtheorem{note}[thm]{Notation}
\newcommand{\norm}[1]{\left\Vert#1\right\Vert}

\begin{document}

\title{Floer homology and surface decompositions}%
\author{Andr\'as Juh\'asz}%
\address{Department of Mathematics, Princeton University, Princeton, NJ 08544, USA}%
\email{ajuhasz@math.princeton.edu}%

\thanks{Research partially supported by OTKA grant no. T49449}
\subjclass{57M27; 57R58}%
\keywords{Sutured manifold; Floer homology; Surface decomposition}

\date{\today}%
%\dedicatory{}%
%\commby{}%
% ----------------------------------------------------------------
\begin{abstract}
Sutured Floer homology, denoted by $SFH,$ is an invariant of
balanced sutured manifolds previously defined by the author. In this
paper we give a formula that shows how this invariant changes under
surface decompositions. In particular, if $(M, \gamma)
\rightsquigarrow (M', \gamma')$ is a sutured manifold decomposition
then $SFH(M',\gamma')$ is a direct summand of $SFH(M, \gamma).$ To
prove the decomposition formula we give an algorithm that computes
$SFH(M,\gamma)$ from a balanced diagram defining $(M,\gamma)$ that
generalizes the algorithm of Sarkar and Wang.

As a corollary we obtain that if $(M, \gamma)$ is taut then
$SFH(M,\gamma) \neq 0.$ Other applications include simple proofs of
a result of Ozsv\'ath and Szab\'o that link Floer homology detects
the Thurston norm, and a theorem of Ni that knot Floer homology
detects fibred knots. Our proofs do not make use of any contact
geometry.

Moreover, using these methods we show that if $K$ is a genus $g$
knot in a rational homology 3-sphere $Y$ whose Alexander polynomial
has leading coefficient $a_g \neq 0$ and if $\text{rk}
\widehat{HFK}(Y,K,g) < 4$ then $Y \setminus N(K)$ admits a depth
$\le 1$ taut foliation transversal to $\partial N(K).$
\end{abstract}

\maketitle
% ----------------------------------------------------------------
\section{Introduction}

In \cite{sutured} we defined a Floer homology invariant for balanced
sutured manifolds. In this paper we study how this invariant changes
under surface decompositions. We need some definitions before we can
state our main result. Recall that $\text{Spin}^c$ structures on
sutured manifolds were defined in \cite{sutured}; all the necessary
definitions can also be found in Section \ref{section:3} of the
present paper.

\begin{defn} \label{defn:1}
Let $(M, \gamma)$ be a balanced sutured manifold and let
$(S,\partial S) \subset (M, \partial M)$ be a properly embedded
oriented surface. An element $\mathfrak{s} \in \text{Spin}^c(M,
\gamma)$ is called \emph{outer} with respect to $S$ if there is a
unit vector field $v$ on $M$ whose homology class is $\mathfrak{s}$
and $v_p \neq -(\nu_S)_p$ for every $p \in S.$ Here $\nu_S$ is the
unit normal vector field of $S$ with respect to some Riemannian
metric on $M$. Let $O_S$ denote the set of outer $\text{Spin}^c$
structures.
\end{defn}

\begin{defn} \label{defn:31}
Suppose that $R$ is a compact, oriented, and open surface. Let $C$
be an oriented simple closed curve in $R.$ If $[C]=0$ in
$H_1(R;\mathbb{Z})$ then $R \setminus C$ can be written as $R_1 \cup
R_2,$ where $R_1$ is the component of $R \setminus C$ that is
disjoint from $\partial R$ and satisfies $\partial R_1 = C.$ We call
$R_1$ the \emph{interior} and $R_2$ the \emph{exterior} of $C.$

We say that the curve $C$ is \emph{boundary-coherent} if either $[C]
\neq 0$ in $H_1(R;\mathbb{Z}),$ or if $[C]=0$ in $H_1(R;\mathbb{Z})$
and $C$ is oriented as the boundary of its interior.
\end{defn}

\begin{thm} \label{thm:1}
Let $(M,\gamma)$ be a balanced sutured manifold and let $(M,
\gamma)\rightsquigarrow^S (M', \gamma')$ be a sutured manifold
decomposition. Suppose that $S$ is open and for every component $V$
of $R(\gamma)$ the set of closed components of $S \cap V$ consists
of parallel oriented boundary-coherent simple closed curves. Then
$$SFH(M', \gamma') = \bigoplus_{\mathfrak{s} \in O_S} SFH(M,\gamma,
\mathfrak{s}).$$ In particular, $SFH(M',\gamma')$ is a direct
summand of $SFH(M,\gamma).$
\end{thm}

In order to prove Theorem \ref{thm:1} we give an algorithm that
computes $SFH(M,\gamma)$ from any given balanced diagram of
$(M,\gamma)$ that generalizes the algorithm of \cite{Sucharit}.

From Theorem \ref{thm:1} we will deduce the following two theorems.
These provide us with positive answers to \cite[Question
9.19]{sutured} and \cite[Conjecture 10.2]{sutured}.

\begin{thm} \label{thm:2}
Suppose that the balanced sutured manifold $(M,\gamma)$ is taut.
Then $$\mathbb{Z} \le SFH(M,\gamma).$$
\end{thm}

If $Y$ is a closed connected oriented 3-manifold and $R \subset Y$
is a compact oriented surface with no closed components then we can
obtain a balanced sutured manifold $Y(R)=(M,\gamma),$ where $M = Y
\setminus \text{Int}(R \times I)$ and $\gamma = \partial R \times
I,$ see \cite[Example 2.6]{sutured}. Furthermore, if $K \subset Y$
is a knot, $\alpha \in H_2(Y,K;\mathbb{Z}),$ and $i \in \mathbb{Z}$
then let $$\widehat{HFK}(Y,K,\alpha,i) = \bigoplus_{\mathfrak{s} \in
\text{Spin}^c(Y,K) \colon \langle c_1(\mathfrak{s}), \alpha \rangle
= 2i} \widehat{HFK}(Y,K,\mathfrak{s}).$$

\begin{thm} \label{thm:3}
Let $K$ be a null-homologous knot in a closed connected oriented
3-manifold $Y$ and let $S \subset Y$ be a Seifert surface of $K.$
Then
$$SFH(Y(S)) \approx \widehat{HFK}(Y,K,[S],g(S)).$$
\end{thm}

\begin{rem} \label{rem:4}
Theorem \ref{thm:3} implies that the invariant $\widehat{HFS}$ of
balanced sutured manifolds defined in \cite{fibred} is equal to
$SFH.$
\end{rem}

Putting these two theorems together we get a new proof of the fact
proved in \cite{OSz6} that knot Floer homology detects the genus of
a knot. In particular, if $Y$ is a rational homology 3-sphere then
$\widehat{HFK}(K,g(K))$ is non-zero and $\widehat{HFK}(K,i) =0$ for
$i > g(K).$

Further applications include a simple proof of a theorem that link
Floer homology detects the Thurston norm, which was proved for links
in $S^3$ in \cite{OSz7}. We generalize this result to links in
arbitrary 3-manifolds. Here we do not use any symplectic or contact
geometry. We also show that the Murasugi sum formula proved in
\cite{Yi} is an easy consequence of Theorem \ref{thm:1}. The main
application of our apparatus is a simplified proof that shows knot
Floer homology detects fibred knots. This theorem was conjectured by
Ozsv\'ath and Szab\'o and first proved in \cite{fibred}. Here we
avoid the contact topology of \cite{Ghiggini} and this allows us to
simplify some of the arguments in \cite{fibred}.

To show the strength of our approach we prove the following
extension of the main result of \cite{fibred}. First we review a few
definitions about foliations, see \cite[Definition 3.8]{Gabai4}.

\begin{defn}
Let $\mathcal{F}$ be a codimension one transversely oriented
foliation. A leaf of $\mathcal{F}$ is of \emph{depth} 0 if it is
compact. Having defined the depth $<p$ leaves we say that a leaf $L$
is depth $p$ if it is proper (i.e., the subspace topology on $L$
equals the leaf topology), $L$ is not of depth $<p,$ and $\bar{L}
\setminus L$ is contained in the union of depth $<p$ leaves. If
$\mathcal{F}$ contains non-proper leaves then the depth of a leaf
may not be defined.

If every leaf of $\mathcal{F}$ is of depth at most $n$ and
$\mathcal{F}$ has a depth $n$ leaf then we say that $\mathcal{F}$ is
\emph{depth} $n.$

A foliation $\mathcal{F}$ is \emph{taut} if there is a single circle
$C$ transverse to $\mathcal{F}$ which intersects every leaf.
\end{defn}

\begin{thm} \label{thm:12}
Let $K$ be a null-homologous genus $g$ knot in a rational homology
3-sphere $Y.$ Suppose that the coefficient $a_g$ of the Alexander
polynomial $\Delta_K(t)$ of $K$ is non-zero and
$$\text{rk} \, \widehat{HFK}(Y,K,g) < 4.$$ Then $Y \setminus
N(K)$ has a depth $\le 1$ taut foliation transverse to $\partial
N(K).$
\end{thm}

\section*{Acknowledgement}

I am grateful for the guidance of Zolt\'an Szab\'o during the course
of this work. I would also like to thank David Gabai, Paolo
Ghiggini, and Yi Ni for the helpful discussions.

\section{Preliminary definitions}

First we briefly review the basic definitions concerning balanced
sutured manifolds and the Floer homology invariant defined for them
in \cite{sutured}.

\begin{defn} \label{defn:2}
A \emph{sutured manifold} $(M,\gamma)$ is a compact oriented
3-manifold $M$ with boundary together with a set $\gamma \subset
\partial M$ of pairwise disjoint annuli $A(\gamma)$ and tori
$T(\gamma).$ Furthermore, the interior of each component of
$A(\gamma)$ contains a \emph{suture}, i.e., a homologically
nontrivial oriented simple closed curve. We denote the union of the
sutures by $s(\gamma).$

Finally every component of $R(\gamma)=\partial M \setminus
\text{Int}(\gamma)$ is oriented. Define $R_+(\gamma)$ (or
$R_-(\gamma)$) to be those components of $\partial M \setminus
\text{Int}(\gamma)$ whose normal vectors point out of (into) $M$.
The orientation on $R(\gamma)$ must be coherent with respect to
$s(\gamma),$ i.e., if $\delta$ is a component of $\partial
R(\gamma)$ and is given the boundary orientation, then $\delta$ must
represent the same homology class in $H_1(\gamma)$ as some suture.
\end{defn}

\begin{defn} \label{defn:3}
A sutured manifold $(M,\gamma)$ is called \emph{balanced} if M has
no closed components, $\chi(R_+(\gamma))=\chi(R_-(\gamma)),$ and the
map $\pi_0(A(\gamma)) \to \pi_0(\partial M)$ is surjective.
\end{defn}

\begin{note} \label{note:3}
Throughout this paper we are going to use the following notation. If
$K$ is a submanifold of the manifold $M$ then $N(K)$ denotes a
regular neighborhood of $K$ in $M.$
\end{note}

For the following see examples 2.3, 2.4, and 2.5 in \cite{sutured}.

\begin{defn} \label{defn:35}
Let $Y$ be a closed connected oriented 3-manifold. Then the balanced
sutured manifold $Y(1)$ is obtained by removing an open ball from
$Y$ and taking an annular suture on its boundary.

Suppose that $L$ is a link in $Y.$ The balanced sutured manifold
$Y(L)=(M,\gamma),$ where $M = Y \setminus N(L)$ and for each
component $L_0$ of $L$ the sutures $\partial N(L_0) \cap s(\gamma)$
consist of two oppositely oriented meridians of $L_0.$

Finally, if $S$ is a Seifert surface in $Y$ then the balanced
sutured manifold $Y(S)=(N,\nu),$ where $N = Y \setminus \text{Int}(S
\times I)$ and $\nu =
\partial S \times I.$
\end{defn}

The following definition can be found for example in
\cite{Scharlemann}.

\begin{defn} \label{defn:36}
Let $S$ be a compact oriented surface (possibly with boundary) whose
components are $S_1, \dots, S_n.$ Then define the \emph{norm} of $S$
to be $$x(S) = \sum_{i \colon \chi(S_i) < 0} |\chi(S_i)|.$$

Let $M$ be a compact oriented 3-manifold and let $N$ be a subsurface
of $\partial M.$ For $s \in H_2(M,N;\mathbb{Z})$ we define its
\emph{norm} $x(s)$ to be the minimum of $x(S)$ taken over all
properly embedded surfaces $(S,\partial S)$ in $(M,N)$ such that
$[S,\partial S] = s.$

If $(S, \partial S) \subset (M,N)$ is a properly embedded oriented
surface then we say that $S$ is \emph{norm minimizing} in $H_2(M,N)$
if $S$ is incompressible and $x(S) = x([S,\partial S])$ for
$[S,\partial S] \in H_2(M,N;\mathbb{Z}).$
\end{defn}

\begin{defn} \label{defn:37}
A sutured manifold $(M,\gamma)$ is \emph{taut} if $M$ is irreducible
and $R(\gamma)$ is norm minimizing in $H_2(M,\gamma).$
\end{defn}

Next we recall the definition of a sutured manifold decomposition,
see \cite[Definition 3.1]{Gabai}.

\begin{defn} \label{defn:4}
Let $(M, \gamma)$ be a sutured manifold. A \emph{decomposing
surface} is a properly embedded oriented surface $S$ in $M$ such
that for every component $\lambda$ of $S \cap \gamma$ one of (1)-(3)
holds:
\begin{enumerate}
\item $\lambda$ is a properly embedded non-separating arc in $\gamma$
such that $|\lambda \cap s(\gamma)| = 1.$
\item $\lambda$ is a simple closed curve in an annular component $A$
of $\gamma$ in the same homology class as $A \cap s(\gamma).$
\item $\lambda$ is a homotopically nontrivial curve in a torus
component $T$ of $\gamma,$ and if $\delta$ is another component of
$T \cap S,$ then $\lambda$ and $\delta$ represent the same homology
class in $H_1(T).$
\end{enumerate}

Then $S$ defines a \emph{sutured manifold decomposition}
$$(M, \gamma)\rightsquigarrow^{S} (M', \gamma'),$$ where $M' = M
\setminus \text{Int}(N(S))$ and $$\gamma' = (\gamma \cap M') \cup
N(S'_+ \cap R_-(\gamma)) \cup N(S'_- \cap R_+(\gamma)), $$
$$R_+(\gamma') = ((R_+(\gamma) \cap M') \cup S'_+) \setminus
\text{Int}(\gamma'),$$
$$R_-(\gamma') = ((R_-(\gamma) \cap M') \cup S'_-) \setminus
\text{Int}(\gamma'),$$ where $S'_+$ ($S'_-$) is the component of
$\partial N(S) \cap M'$ whose normal vector points out of (into)
$M'.$
\end{defn}

\begin{defn} \label{defn:33}
A decomposing surface $S$ in $(M,\gamma)$ is called a \emph{product
disk} if $S$ is a disk such that $|D \cap s(\gamma)| = 2.$ A surface
decomposition $(M,\gamma) \rightsquigarrow^S (M',\gamma')$ is called
a \emph{product decomposition} if $S$ is a product disk.
\end{defn}

\begin{defn} \label{defn:34}
A decomposing surface $S$ lying in the sutured manifold $(M,\gamma)$
is called a \emph{product annulus} if $S$ is an annulus, one
component of $\partial S$ is contained in $R_+(\gamma),$ and the
other component is contained in $R_-(\gamma).$
\end{defn}

\begin{defn} \label{defn:5}
A \emph{sutured Heegaard diagram} is a tuple $( \Sigma,
\boldsymbol{\alpha}, \boldsymbol{\beta}),$ where $\Sigma$ is a
compact oriented surface with boundary and $\boldsymbol{\alpha}$ and
$\boldsymbol{\beta}$ are two sets of pairwise disjoint simple closed
curves in $\text{Int}(\Sigma).$
\end{defn}

Every sutured Heegaard diagram $( \Sigma, \boldsymbol{\alpha},
\boldsymbol{\beta})$ uniquely \emph{defines} a sutured manifold $(M,
\gamma)$ using the following construction. Suppose that
$\boldsymbol{\alpha}=\{\,\alpha_1,\dots,\alpha_m\,\}$ and
$\boldsymbol{\beta}=\{\,\beta_1,\dots,\beta_n\,\}.$ Let $M$ be the
3-manifold obtained from $\Sigma \times I$ by attaching
3-dimensional 2-handles along the curves $\alpha_i \times \{0\}$ and
$\beta_j \times \{1\}$ for $i=1, \dots, m$ and $j=1, \dots, n.$ The
sutures are defined by taking $\gamma =
\partial \Sigma \times I$ and $s(\gamma)=
\partial \Sigma \times \{1/2\}.$

\begin{defn} \label{defn:6}
A sutured Heegaard diagram $( \Sigma, \boldsymbol{\alpha},
\boldsymbol{\beta})$ is called \emph{balanced} if
$|\boldsymbol{\alpha}| = |\boldsymbol{\beta}|$ and the maps
$\pi_0(\partial \Sigma) \to \pi_0(\Sigma \setminus \bigcup
\boldsymbol{\alpha})$ and $\pi_0(\partial \Sigma) \to \pi_0(\Sigma
\setminus \bigcup \boldsymbol{\beta})$ are surjective.
\end{defn}

The following is \cite[Proposition 2.14]{sutured}.

\begin{prop} \label{prop:1}
For every balanced sutured manifold $(M,\gamma)$ there exists a
balanced diagram defining it.
\end{prop}

\begin{defn} \label{defn:7}
For a balanced diagram let $\mathcal{D}_1, \dots, \mathcal{D}_m$
denote the closures of the components of $\Sigma \setminus (\bigcup
\boldsymbol{\alpha} \cup \bigcup \boldsymbol{\beta})$ disjoint from
$\partial \Sigma.$ Then let $D(\Sigma, \boldsymbol{\alpha},
\boldsymbol{\beta})$ be the free abelian group generated by
$\{\,\mathcal{D}_1, \dots, \mathcal{D}_m\,\}.$ This is of course
isomorphic to $\mathbb{Z}^m.$ We call an element of $D(\Sigma,
\boldsymbol{\alpha}, \boldsymbol{\beta})$ a \emph{domain}. An
element $\mathcal{D}$ of $\mathbb{Z}_{\ge 0}^m$ is called a
\emph{positive} domain, we write $\mathcal{D} \ge 0.$ A domain
$\mathcal{P} \in D(\Sigma, \boldsymbol{\alpha}, \boldsymbol{\beta})$
is called a \emph{periodic domain} if the boundary of the 2-chain
$\mathcal{P}$ is a linear combination of full $\alpha$- and
$\beta$-curves.
\end{defn}

\begin{defn} \label{defn:8}
A balanced diagram $(\Sigma, \boldsymbol{\alpha},
\boldsymbol{\beta})$ is called \emph{admissible} if every periodic
domain $\mathcal{P} \neq 0$ has both positive and negative
coefficients.
\end{defn}

The following proposition is \cite[Corollary 3.12]{sutured}.

\begin{prop} \label{prop:2}
If $(M, \gamma)$ is a balanced sutured manifold such that $$H_2(M;
\mathbb{Z}) = 0$$ and if $(\Sigma, \boldsymbol{\alpha},
\boldsymbol{\beta})$ is an arbitrary balanced diagram defining
$(M,\gamma)$ then there are no non-zero periodic domains in
$D(\Sigma, \boldsymbol{\alpha, \boldsymbol{\beta}})$. Thus any
balanced diagram defining $(M, \gamma)$ is automatically admissible.
\end{prop}

For a surface $\Sigma$ let $\text{Sym}^d(\Sigma)$ denote the d-fold
symmetric product $\Sigma^{\times d} / S_d.$ This is a smooth
$2d$-manifold. A complex structure $\mathfrak{j}$ on $\Sigma$
naturally endows $\text{Sym}^d(\Sigma)$ with a complex structure.
Let $(\Sigma, \boldsymbol{\alpha}, \boldsymbol{\beta})$ be a
balanced diagram, where $\boldsymbol{\alpha} = \{\, \alpha_1, \dots,
\alpha_d \,\}$ and $\boldsymbol{\beta} = \{\, \beta_1, \dots,
\beta_d \,\}.$ Then the tori $\mathbb{T}_{\alpha} = (\alpha_1 \times
\dots \times \alpha_d) / S_d$ and $\mathbb{T}_{\beta} = (\beta_1
\times \dots \times \beta_d) / S_d$ are $d$-dimensional totally real
submanifolds of $\text{Sym}^d(\Sigma).$

\begin{defn} \label{defn:9}
Let $\mathbf{x}, \mathbf{y} \in \mathbb{T}_{\alpha} \cap
\mathbb{T}_{\beta}.$ A domain $\mathcal{D} \in D(\Sigma,
\boldsymbol{\alpha}, \boldsymbol{\beta})$ is said to \emph{connect
$\mathbf{x}$ to $\mathbf{y}$} if for every $1 \le i \le d$ the
equalities $\partial(\alpha_i \cap \partial \mathcal{D}) =
(\mathbf{x} \cap \alpha_i) - (\mathbf{y} \cap \alpha_i)$ and
$\partial(\beta_i \cap \partial \mathcal{D}) = (\mathbf{x} \cap
\beta_i) - (\mathbf{y} \cap \beta_i)$ hold. We are going to denote
by $D(\mathbf{x}, \mathbf{y})$ the set of domains connecting
$\mathbf{x}$ to $\mathbf{y}.$
\end{defn}

\begin{note} \label{note:1}
Let $\mathbb{D}$ denote the unit disc in $\mathbb{C}$ and let
$e_1=\{\, z \in \partial \mathbb{D} \colon \text{Re}(z) \ge 0\,\}$
and $e_2=\{\, z \in \partial \mathbb{D} \colon \text{Re}(z) \le
0\,\}.$
\end{note}

\begin{defn} \label{defn:10}
Let $\mathbf{x}, \mathbf{y} \in \mathbb{T}_{\alpha} \cap
\mathbb{T}_{\beta}$ be intersection points. A \emph{Whitney disc
connecting $\mathbf{x}$ to $\mathbf{y}$} is a continuous map $u
\colon \mathbb{D} \to \text{Sym}^d(\Sigma)$ such that $u(-i) =
\mathbf{x},$ $u(i) = \mathbf{y}$ and $u(e_1) \subset
\mathbb{T}_{\alpha},$ $u(e_2) \subset \mathbb{T}_{\beta}.$ Let
$\pi_2(\mathbf{x}, \mathbf{y})$ denote the set of homotopy classes
of Whitney discs connecting $\mathbf{x}$ to $\mathbf{y}.$
\end{defn}

\begin{defn}
If $z \in \Sigma \setminus (\bigcup \boldsymbol{\alpha} \cup \bigcup
\boldsymbol{\beta})$ and if $u$ is a Whitney disc then choose a
Whitney disc $u'$ homotopic to $u$ such that $u'$ intersects the
hypersurface $ \{z\} \times \text{Sym}^{d-1}(\Sigma)$ transversally.
Define $n_z(u)$ to be the algebraic intersection number $u' \cap
(\{z\} \times \text{Sym}^{d-1}(\Sigma)).$
\end{defn}

\begin{defn}
Let $\mathcal{D}_1, \dots, \mathcal{D}_m$ be as in Definition
\ref{defn:7}. For every $1 \le i \le m$ choose a point $z_i \in
\mathcal{D}_i.$ Then the \emph{domain of a Whitney disc $u$} is
defined as
$$\mathcal{D}(u) = \sum_{i=1}^m n_{z_i}(u)\mathcal{D}_i \in
D(\Sigma, \boldsymbol{\alpha}, \boldsymbol{\beta}).$$ If $\phi \in
\pi_2(\mathbf{x}, \mathbf{y})$ and if $u$ is a representative of the
homotopy class $\phi$ then let $\mathcal{D}(\phi) = \mathcal{D}(u).$
\end{defn}

\begin{defn} \label{defn:11}
We define the Maslov index of a domain $\mathcal{D} \in D(\Sigma,
\boldsymbol{\alpha}, \boldsymbol{\beta})$ as follows. If there is a
homotopy class $\phi$ of Whitney discs such that $\mathcal{D}(\phi)
= \mathcal{D}$ then let $\mu(\mathcal{D}) = \mu(\phi).$ Otherwise we
define $\mu(\mathcal{D})$ to be $-\infty.$ Furthermore, let
$\mathcal{M}(\mathcal{D})$ denote the moduli space of holomorphic
Whitney discs $u$ such that $\mathcal{D}(u) = \mathcal{D}$ and let
$\widehat{\mathcal{M}}(\mathcal{D}) =
\mathcal{M}(\mathcal{D})/\mathbb{R}.$
\end{defn}

Let $(M, \gamma)$ be a balanced sutured manifold and $(\Sigma,
\boldsymbol{\alpha}, \boldsymbol{\beta})$ an admissible balanced
diagram defining it. Fix a coherent system of orientations as in
\cite[Definition 3.11]{OSz}. Then for a generic almost complex
structure each moduli space $\widehat{\mathcal{M}}(\mathcal{D})$ is
a compact oriented manifold of dimension $\mu(\mathcal{D})-1.$ We
denote by $CF(\Sigma, \boldsymbol{\alpha}, \boldsymbol{\beta})$ the
free abelian group generated by the points of $\mathbb{T}_{\alpha}
\cap \mathbb{T}_{\beta}.$ We define an endomorphism $\partial \colon
CF(\Sigma, \boldsymbol{\alpha}, \boldsymbol{\beta}) \to CF(\Sigma,
\boldsymbol{\alpha}, \boldsymbol{\beta})$ such that on each
generator $\mathbf{x} \in \mathbb{T}_{\alpha} \cap
\mathbb{T}_{\beta}$ it is given by the formula
$$\partial \mathbf{x} = \sum_{\mathbf{y} \in \mathbb{T}_{\alpha}
\cap \,\mathbb{T}_{\beta}} \sum_{\{\,\mathcal{D} \in D(\mathbf{x},
\mathbf{y})\colon \mu(\mathcal{D})=1\,\}} \#
\widehat{\mathcal{M}}(\mathcal{D})\cdot \mathbf{y}.$$ Then
$(CF(\Sigma, \boldsymbol{\alpha}, \boldsymbol{\beta}),\partial)$ is
a chain complex whose homology depends only on the underlying
sutured manifold $(M,\gamma).$ We denote this homology group by
$SFH(M,\gamma).$

For the following see \cite[Proposition 9.1]{sutured} and
\cite[Proposition 9.2]{sutured}.

\begin{prop}
If $Y$ is a closed connected oriented 3-manifold then
$$SFH(Y(1)) \approx \widehat{HF}(Y).$$ Furthermore, if $L$ is a link
in $Y$ and $\vec{L}$ is an arbitrary orientation of $L$ then
$$SFH(Y(L)) \otimes \mathbb{Z}_2 \approx \widehat{HFL}(\vec{L}).$$
\end{prop}

\section{$\text{Spin}^c$ structures and relative Chern classes}
\label{section:3}

First we review the definition of a $\text{Spin}^c$ structure on a
balanced sutured manifold $(M,\gamma)$ that was introduced in
\cite{sutured}. Note that in a balanced sutured manifold none of the
sutures are tori. Fix a Riemannian metric on $M.$

\begin{note} \label{note:2}
Let $v_0$ be a nowhere vanishing vector field along $\partial M$
that points into $M$ along $R_-(\gamma),$ points out of $M$ along
$R_+(\gamma),$ and on $\gamma$ it is the gradient of the height
function $s(\gamma) \times I \to I.$ The space of such vector fields
is contractible.
\end{note}

\begin{defn} \label{defn:12}
Let $v$ and $w$ be nowhere vanishing vector fields on $M$ that agree
with $v_0$ on $\partial M.$ We say that $v$ and $w$ are
\emph{homologous} if there is an open ball $B \subset \text{Int}(M)$
such that $v|(M \setminus B)$ is homotopic to $w|(M \setminus B)$
through nowhere vanishing vector fields rel $\partial M.$ We define
$\text{Spin}^c(M, \gamma)$ to be the set of homology classes of
nowhere vanishing vector fields $v$ on $M$ such that $v|\partial M =
v_0.$
\end{defn}

\begin{defn} \label{defn:20}
Let $(M, \gamma)$ be a balanced sutured manifold and $(\Sigma,
\boldsymbol{\alpha},\boldsymbol{\beta})$ a balanced diagram defining
it. To each $\mathbf{x} \in \mathbb{T}_{\alpha} \cap
\mathbb{T}_{\beta}$ we assign a $\text{Spin}^c$ structure
$\mathfrak{s}(\mathbf{x}) \in \text{Spin}^c(M, \gamma)$ as follows.
Choose a Morse function $f$ on $M$ compatible with the given
balanced diagram $(\Sigma, \boldsymbol{\alpha},
\boldsymbol{\beta}).$ Then $\mathbf{x}$ corresponds to a
multi-trajectory $\gamma_{\mathbf{x}}$ of $\text{grad}(f)$
connecting the index one and two critical points of $f$. In a
regular neighborhood $N(\gamma_{\mathbf{x}})$ we can modify
$\text{grad}(f)$ to obtain a nowhere vanishing vector field $v$ on
$M$ such that $v|\partial M = v_0.$ We define
$\mathfrak{s}(\mathbf{x})$ to be the homology class of this vector
field $v.$
\end{defn}

\begin{prop} \label{prop:3}
The vector bundle $v_0^{\perp}$ over $\partial M$ is trivial if and
only if for every component $F$ of $\partial M$ the equality $\chi(F
\cap R_+(\gamma)) = \chi(F \cap R_-(\gamma))$ holds.
\end{prop}

\begin{proof}
Since $v_0^{\perp}|R_+(\gamma) = TR_+(\gamma)$ and
$v_0^{\perp}|R_-(\gamma) = -TR_-(\gamma)$ we get that
$$\left\langle\, e(v_0^{\perp}|F), [F] \,\right\rangle = \chi(F \cap
R_+(\gamma)) - \chi(F \cap R_-(\gamma)).$$ Furthermore, the rank two
bundle $v_0^{\perp}|F$ is trivial if and only if its Euler class
vanishes.
\end{proof}

\begin{defn}
We call a sutured manifold $(M,\gamma)$ \emph{strongly balanced} if
for every component $F$ of $\partial M$ the equality $\chi(F \cap
R_+(\gamma)) = \chi(F \cap R_-(\gamma))$ holds.
\end{defn}

\begin{rem} \label{rem:1}
Note that if $(M,\gamma)$ is balanced then we can associate to it a
strongly balanced sutured manifold $(M',\gamma')$ such that
$(M,\gamma)$ can be obtained from $(M', \gamma')$ by a sequence of
product decompositions. We can construct such an $(M', \gamma')$ as
follows. If $F_1$ and $F_2$ are distinct components of $\partial M$
then choose two points $p_1 \in s(\gamma) \cap F_1$ and $p_2 \in
s(\gamma) \cap F_2.$ For $i =1,2$ let $D_i$ be a small neighborhood
of $p_i$ homeomorphic to a closed disc. We get a new sutured
manifold by gluing together $D_1$ and $D_2.$ Then $(M,\gamma)$ can
be retrieved by decomposing along $D_1 \sim D_2.$ By repeating this
process we get a sutured manifold $(M', \gamma')$ with a single
boundary component. Since $(M,\gamma)$ was balanced $(M', \gamma')$
is strongly balanced. By adding such product one-handles we can even
achieve that $\gamma$ is connected.
\end{rem}

\begin{defn} \label{defn:13}
Suppose that $(M,\gamma)$ is a strongly balanced sutured manifold.
Let $t$ be a trivialization of $v_0^{\perp}$ and let $\mathfrak{s}
\in \text{Spin}^c(M,\gamma).$ Then we define $$c_1(\mathfrak{s},t)
\in H^2(M, \partial M; \mathbb{Z})$$ to be the relative Euler class
of the vector bundle $v^{\perp}$ with respect to the trivialization
$t.$ In other words, $c_1(\mathfrak{s},t)$ is the obstruction to
extending $t$ from $\partial M$ to a trivialization of $v^{\perp}$
over $M.$
\end{defn}

\begin{defn} \label{defn:14}
Let $S$ be a decomposing surface in a balanced sutured manifold
$(M,\gamma)$ such that the positive unit normal field $\nu_S$ of $S$
is nowhere parallel to $v_0$ along $\partial S.$ This holds for
generic $S.$ We endow $\partial S$ with the boundary orientation.
Let us denote the components of $\partial S$ by $T_1, \dots, T_k.$

Let $w_0$ denote the projection of $v_0$ into $TS,$ this is a
nowhere zero vector field. Moreover, let $f$ be the positive unit
tangent vector field of $\partial S.$ For $1 \le i \le k$ we define
the \emph{index} $I(T_i)$ to be the number of times $w_0$ rotates
with respect to $f$ as we go around $T_i.$ Then define $$I(S) =
\sum_{i=1}^k I(T_k).$$

Let $p(\nu_S)$ be the projection of $\nu_S$ into $v^{\perp}.$
Observe that $p(\nu_S)|\partial S$ is nowhere zero. For $1 \le i \le
k$ we define $r(T_i,t)$ to be the rotation of $p(\nu_S)|\partial
T_i$ with respect to the trivialization $t$ as we go around $T_i.$
Moreover, let $$r(S,t) = \sum_{i=1}^k r(T_i,t).$$

We introduce the notation
$$c(S,t) = \chi(S) + I(S) - r(S,t).$$
\end{defn}

\begin{lem} \label{lem:9}
Let $(M,\gamma)$ be a balanced sutured manifold and let $S$ be a
decomposing surface as in Definition \ref{defn:14}.
\begin{enumerate}
\item If $T$ is a component of $\partial S$ such that $T \not\subset
\gamma$ then $$I(T) = -\frac{|T \cap s(\gamma)|}{2}.$$
\item Suppose that $T_1, \dots, T_a$ are components of $\partial S$
such that $\mathcal{T} = T_1 \cup \dots \cup T_a \subset \gamma$ is
parallel to $s(\gamma)$ and $\nu_S$ points out of $M$ along
$\mathcal{T}.$ Then $I(T_j) = 0$ for $1 \le j \le a;$ moreover,
$$\sum_{j=1}^a r(T_j,t) = \chi(R_+(\gamma)).$$
\end{enumerate}
\end{lem}

\begin{proof}
First we prove part (1). We can suppose that $w_0$ is tangent to $T$
exactly at the points of $\partial T \cap s(\gamma).$ Then at a
point $p \in T \cap s(\gamma)$ we have $w_0/|w_0| = f$ if and only
if $T$ goes from $R_-(\gamma)$ to $R_+(\gamma)$ and in that case
$w_0$ rotates from the inside of $S$ to the outside, see Figure
\ref{fig:7}. Thus $w_0$ rotates $-|T \cap s(\gamma)|/2$ times with
respect to $f$ as we go around $T.$

\begin{figure}[t]
\includegraphics{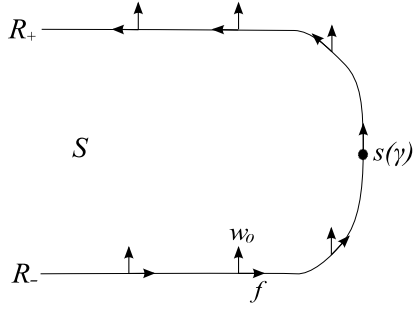}
\caption{If $T \not\subset \gamma$ then the index $I(T)$ is $-|T
\cap s(\gamma)|/2.$} \label{fig:7}
\end{figure}

Now we prove part (2). Let $1 \le j \le a.$ Since $\nu_S$ points out
of $M$ along $T_j$ we get that $w_0$ points into $S$ along $T_j.$ So
$w_0$ and $f$ are nowhere equal along $T_j,$ and thus $I(T_j)=0.$

Since $\mathcal{T}$ is parallel to $s(\gamma)$ it bounds a surface
$\mathcal{R}_+ \subset \partial M$ which is diffeomorphic to
$R_+(\gamma)$ and contains $R_+(\gamma).$ Since $\nu_S$ points out
of $M$ along $\mathcal{T}$ there is an isomorphism $i \colon
v_0^{\perp}|\mathcal{R}_+ \to T\mathcal{R}_+$ such that
$i(p(\nu_S))$ is an outward normal field of $\mathcal{R}_+$ along
$\partial \mathcal{R}_+.$ Moreover, $i(t|\mathcal{R}_+)$ gives a
trivialization of $T\mathcal{R}_+.$ Using the Poincar\'e-Hopf
theorem we get that $p(\nu_S)$ rotates $\chi(\mathcal{R}_+) =
\chi(R_+(\gamma))$ times with respect to $t$ as we go around
$\mathcal{T}.$
\end{proof}

Recall that we defined the notion of an outer $\text{Spin}^c$
structure in Definition \ref{defn:1}.

\begin{lem} \label{lem:1}
Suppose that $(M,\gamma)$ is a strongly balanced sutured manifold.
Let $t$ be a trivialization of $v_0^{\perp}$, let $\mathfrak{s} \in
\text{Spin}^c(M,\gamma),$ and let $S$ be a decomposing surface in
$(M, \gamma)$ as in Definition \ref{defn:14}. Then $\mathfrak{s}$ is
outer with respect to $S$ if and only if
\begin{equation} \label{eqn:1}
\left\langle\, c_1(\mathfrak{s},t), [S] \,\right\rangle = c(S,t).
\end{equation}
\end{lem}

\begin{proof}
Endow $M$ with an arbitrary Riemannian metric. First we show that if
$\mathfrak{s} \in O_S$ then equation \ref{eqn:1} holds. Using the
naturality of Chern classes it is sufficient to prove that if $v$ is
a unit vector field over $S$ that agrees with $v_0$ over $\partial
S$ and is nowhere equal to $-\nu_S$ then $\langle\,
c_1(v^{\perp},t), [S] \,\rangle = c(S,t).$

If we project $\nu_S$ into $v^{\perp}$ we get a section $p(\nu_S)$
of $v^{\perp}$ that vanishes exactly where $\nu_S = v.$ We can
perturb $v$ slightly to make all tangencies between $v^{\perp}$ and
$S$ non-degenerate. Let $e$ and $h$ denote the number of elliptic,
respectively hyperbolic tangencies between $v^{\perp}$ and $S.$ At
each such tangency the orientation of $v^{\perp}$ and $TS$ agree.
Thus $\langle\, c_1(v^{\perp},t_1), [S] \,\rangle = e-h,$ where $t_1
= p(\nu_S)|\partial S.$ Since
$$\left\langle\, c_1(v^{\perp},t_1) - c_1(v^{\perp},t),[S]
\,\right\rangle = r(S,t)$$ we get that
$$\left\langle\, c_1(v^{\perp},t), [S] \,\right\rangle = e-h - r(S,t).$$

On the other hand, if we project $v$ into $TS$ we get a vector field
$w$ on $S$ that is zero exactly at the points where $\nu_S = v$ as
well. Note that $w$ has index $1$ exactly where $v^{\perp}$ and $S$
have an elliptic tangency and has index $-1$ at hyperbolic
tangencies. Moreover, $w|\partial S = w_0.$ If we extend $f$ to a
vector field $f_1$ over $S$ the sum of the indices of $f_1$ will by
$\chi(S)$ by the Poincar\'e-Hopf theorem. Putting these observations
together we get that
$$I(S) = (e-h) - \chi(S).$$ So we conclude that
$$\left\langle\,c_1(v^{\perp},t), [S] \,\right\rangle = \chi(S)
+ I(S) - r(S,t) = c(S,t).$$

Now we prove that if for $\mathfrak{s} \in \text{Spin}^c(M,\gamma)$
equation \ref{eqn:1} holds then $\mathfrak{s} \in O_S.$ Let $STM$
denote the unit sphere bundle of $TM.$ Then $v_0|\partial S$ is a
section over $\partial S$ of $(STM|S) \setminus (-\nu_S),$ which is
a bundle over $S$ with contractible fibers. Thus $v_0|\partial S$
extends to a section $v_1 \colon S \to STM|S$ that is nowhere equal
to $-\nu_S.$ In the first part of the proof we showed that for such
a vector field $v_1$ the equation $\langle\, c_1(v_1^{\perp},t), [S]
\,\rangle = c(S,t)$ holds.

Let $v'$ be a unit vector field over $M$ whose homology class is
$\mathfrak{s}$ and let $v=v'|S.$ Since $\mathfrak{s}$ satisfies
equation \ref{eqn:1} we get that $$\left\langle\, c_1(v^{\perp},t)
-c_1(v_1^{\perp},t), [S] \,\right\rangle = 0.$$ The obstruction
class $o(v, v_1) \in H^2(S, \partial S; \mathbb{Z})$ vanishes if and
only if the sections $v$ and $v_1$ of $STM|S$ are homotopic relative
to $\partial S.$ A cochain $o$ representing $o(v,v_1)$ can be
obtained as follows. First take a triangulation of $S$ and a
trivialization of $STM|S.$ Then $v$ and $v_1$ can be considered to
be maps from $S$ to $S^2.$ One can homotope $v$ rel $\partial S$ to
agree with $v_1$ on the one-skeleton of $S.$ The value of $o$ on a
2-simplex $\Delta$ is the difference of $v|\Delta$ and $v_1|\Delta,$
which is an element of $\pi_2(S^2) \approx \mathbb{Z}.$ Since
$2o(v,v_1) = c_1(v^{\perp},t) -c_1(v_1^{\perp},t)$ and $H^2(S,
\partial S; \mathbb{Z})$ is torsion free we get that $o(v,v_1) =0,$
i.e., $v$ is homotopic to $v_1$ rel $\partial S.$ By extending this
homotopy of $v'$ fixing $v'|\partial M$ we get a vector field $v_1'$
on $M$ that agrees with $v_1$ on $S.$ Thus $\mathfrak{s}$ can be
represented by the vector field $v_1'$ that is nowhere equal to
$-\nu_S,$ and so $\mathfrak{s} \in O_S.$
\end{proof}

In light of Lemma \ref{lem:1} we can reformulate Theorem \ref{thm:1}
for strongly balanced sutured manifolds as follows.

\begin{thm} \label{thm:4}
Let $(M,\gamma)$ be a strongly balanced sutured manifold;
furthermore, let $(M, \gamma)\rightsquigarrow^S (M', \gamma')$ be a
sutured manifold decomposition. Suppose that $S$ is open and for
every component $V$ of $R(\gamma)$ the set of closed components of
$S \cap V$ consists of parallel oriented boundary-coherent simple
closed curves. Choose a trivialization $t$ of $v_0^{\perp}.$ Then
$$SFH(M', \gamma') = \bigoplus_{\mathfrak{s} \in
\text{Spin}^c(M,\gamma) \colon \langle\,
c_1(\mathfrak{s},t),[S]\,\rangle = c(S,t)} SFH(M,\gamma,
\mathfrak{s}).$$
\end{thm}

\section{Finding a balanced diagram adapted to a decomposing surface}

\begin{defn} \label{defn:15}
We say that the decomposing surfaces $S_0$ and $S_1$ are
\emph{equivalent} if they can be connected by an isotopy through
decomposing surfaces.
\end{defn}

\begin{rem} \label{rem:2}
During an isotopy through decomposing surfaces the number of arcs of
$S \cap \gamma$ can never change. Moreover, if $S_0$ and $S_1$ are
equivalent then decomposing along them give the same sutured
manifold.
\end{rem}

\begin{defn} \label{defn:16}
A balanced diagram \emph{adapted} to the decomposing surface $S$ in
$(M,\gamma)$ is a quadruple $(\Sigma, \boldsymbol{\alpha},
\boldsymbol{\beta}, P)$ that satisfies the following conditions.
$(\Sigma, \boldsymbol{\alpha}, \boldsymbol{\beta})$ is a balanced
diagram of $(M,\gamma);$ furthermore, $P \subset \Sigma$ is a
quasi-polygon (i.e., a closed subsurface of $\Sigma$ with polygonal
boundary) such that $P \cap \partial \Sigma$ is exactly the set of
vertices of $P.$ We are also given a decomposition $\partial P = A
\cup B,$ where both $A$ and $B$ are unions of pairwise disjoint
edges of $P.$ This decomposition has to satisfy the property that
$\alpha \cap B = \emptyset$ and $\beta \cap A = \emptyset$ for every
$\alpha \in \boldsymbol{\alpha}$ and $\beta \in \boldsymbol{\beta}.$
Finally, $S$ is given up to equivalence by smoothing the corners of
the surface $(P\times \{1/2\}) \cup (A \times [1/2,1]) \cup (B
\times [0,1/2]) \subset (M, \gamma)$ (see Definition \ref{defn:5}).
The orientation of $S$ is given by the orientation of $P \subset
\Sigma.$ We call a tuple $(\Sigma, \boldsymbol{\alpha},
\boldsymbol{\beta}, P)$ satisfying the above conditions a
\emph{surface diagram}.
\end{defn}

\begin{prop} \label{prop:4}
Suppose that $S$ is a decomposing surface in the balanced sutured
manifold $(M,\gamma).$ If the boundary of each component of $S$
intersects both $R_+(\gamma)$ and $R_-(\gamma)$ (in particular $S$
is open) and $\partial S$ has no closed component lying entirely in
$\gamma$ then there exists a Heegaard diagram of $(M,\gamma)$
adapted to $S.$
\end{prop}

\begin{proof}
We are going to construct a self-indexing Morse function $f$ on $M$
with no minima and maxima as in the proof of \cite[Proposition
2.13]{sutured} with some additional properties. In particular, we
require that $f|R_-(\gamma) \equiv -1$ and $f|R_+(\gamma)\equiv 4.$
Furthermore, $f|\gamma$ is given by the formula $p_2 \circ \varphi,$
where $\varphi \colon \gamma \to s(\gamma) \times [-1,4]$ is a
diffeomorphism such that $\varphi(s(\gamma)) = s(\gamma) \times
\{3/2\}$ and $p_2 \colon s(\gamma) \times [-1,4] \to [-1,4]$ is the
projection onto the second factor. We choose $\varphi$ such that
each arc of $S \cap \gamma$ maps to a single point under $p_1 \circ
\varphi \colon \gamma \to s(\gamma).$

\begin{figure}[t]
\includegraphics{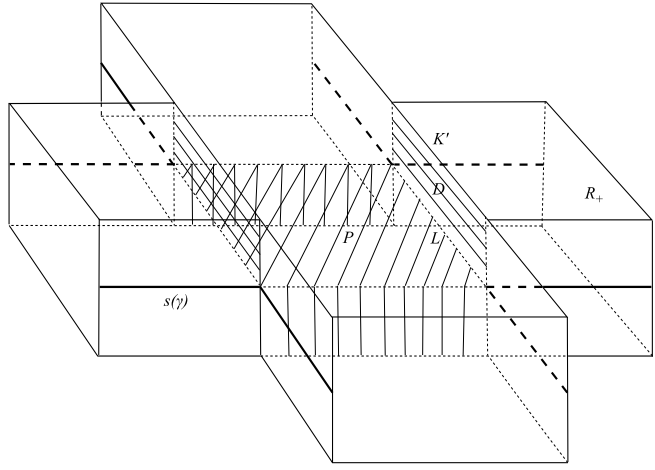}
\caption{This diagram shows a decomposing surface which is a disk
that intersects $s(\gamma)$ in four points.}  \label{fig:1}
\end{figure}

We are going to define a quasi-polygon $P \subset S$ such that $S
\cap s(\gamma)$ is the set of vertices of $P,$ see Figure
\ref{fig:1}. Let $K_1, \dots, K_{m+n}$ be the closures of the
components of $\partial S \setminus s(\gamma)$ enumerated such that
$K_i$ is an arc for $1 \le i \le m$ and $K_i$ is a circle for $m+1
\le i \le m+n.$

For every $1 \le i \le m$ choose an arc $L_i$ whose interior lies in
$\text{int}(S)$ parallel to $K_i$ and such that $\partial L_i
=\partial K_i.$ Moreover, let $D_i$ be the closed bigon bounded by
$K_i$ and $L_i$ and define $K_i' = K_i \cap R(\gamma).$ Also choose
a diffeomorphism $d_i \colon D_i \to I \times I$ that takes $K_i'$
to $I \times \{0\}$ and $L_i$ to $I \times \{1\}$ and such that for
each $t \in [0,1]$ we have $f \circ d_i^{-1}(0,t) = f \circ
d_i^{-1}(1,t).$ Note that $f$ is already defined on $\partial M.$ We
define $f$ on $D_i$ by the formula
$$f(d_i^{-1}(u,t)) = f(d_i^{-1}(0,t)).$$

If $m+1 \le i \le m+n$ then let $L_i$ be a circle parallel to $K_i$
lying in the interior of $S.$ Let $D_i$ be the annulus bounded by
$K_i$ and $L_i.$ Choose a diffeomorphism $$d_i \colon D_i \to S^1
\times J_i,$$ where $J_i = [3/2,4]$ if $K_i \subset R_+(\gamma)$ and
$J_i = [-1,3/2]$ otherwise. In both cases we require that $d_i(L_i)
= 3/2.$ Then let $f|D_i = \pi_2 \circ d_i,$ where $\pi_2 \colon S^1
\times J_i \to J_i$ is the projection onto the second factor.

We take $$\partial P = \bigcup_{i=1}^{m+n} L_i,$$ and $L_i$ will be
an edge of $\partial P$ for every $1 \le i \le m+n.$ The
decomposition $\partial P = A \cup B$ is given by taking $A$ to be
the union of those edges $L_i$ of $\partial P$ for which $K_i \cap
R_+(\gamma) \neq \emptyset.$

Let $P$ be the closure of the component of $S \setminus \partial P$
that is disjoint from $\partial S.$ For $p \in P$ let $f(p) = 3/2.$
Note that the function $f|S$ is not smooth along $\partial P,$ so we
modify $S$ by introducing a right angle edge along $\partial P$
(such that we get back $S$ after smoothing the corners). There are
essentially two ways of creasing $S$ along an edge $L_i$ of $P.$ Let
$\nu_P = \nu_S|P$ be the positive unit normal field of $P$ in $M.$
If $L_i \subset A$ then we choose the crease such that $\nu_P|L_i$
points into $D_i$ and if $L_i \subset B$ then we require that
$\nu_P|L_i$ points out of $D_i.$

Now extend $f$ from $\partial M \cup S$ to a Morse function $f_0$ on
$M.$ Then $$P = S \cap f_0^{-1}(3/2).$$ We choose the extension
$f_0$ as follows. For $1 \le i \le m+n$ let $N(D_i)$ be a regular
neighborhood of $D_i$ and let $T_i \colon N(D_i) \to D_i \times
[-1,1]$ be a diffeomorphism. Then for $(x,t) \in D_i \times [-1,1]$
let
$$f_0(T_i^{-1}(x,t)) = f(x).$$ Due to the choice of the creases we can
define $f_0$ such that $\text{grad}(f)|P \neq -\nu_S.$ Thus we have
achieved that for each $a \in A$ the gradient flow line of $f_0$
coming out of $a$ ends on $R_+(\gamma)$ and for each $b \in B$ the
negative gradient flow line of $f_0$ going through $b$ ends on
$R_-(\gamma).$

By making $f_0$ self-indexing we obtain a Morse function $f.$
Suppose that the Heegaard diagram corresponding to $f$ is $(\Sigma,
\boldsymbol{\alpha}, \boldsymbol{\beta}).$ We have two partitions
$\boldsymbol{\alpha} = \boldsymbol{\alpha}_0 \cup
\boldsymbol{\alpha}_1$ and $\boldsymbol{\beta} =
\boldsymbol{\beta}_0 \cup \boldsymbol{\beta}_1,$ where curves in
$\boldsymbol{\alpha}_1$ correspond to index one critical points $p$
of $f_0$ for which $f_0(p) > 3/2$ and $\boldsymbol{\beta}_1$ comes
from those index two critical points $q$ of $f_0$ for which $f_0(q)
< 3/2.$ Then $f^{-1}(3/2)$ differs from $f_0^{-1}(3/2)$ as follows.
Add an $S^2$ component to $f_0^{-1}(3/2)$ for each index zero
critical point of $f_0$ lying above $3/2$ and for each index three
critical point of $f_0$ lying below $3/2.$ Then add two-dimensional
one-handles to the previous surface whose belt circles are the
curves in $\boldsymbol{\alpha}_1 \cup \boldsymbol{\beta}_1.$
%Then $f^{-1}(3/2)$ differs from $f_0^{-1}(3/2)$ by adding
%some $S^2$ components corresponding to index zero critical points of
%$f_0$ lying above $3/2$ and index three critical points of $f_0$
%lying below $3/2,$ plus by the addition of a two-dimensional
%one-handle for each curve in $\boldsymbol{\alpha}_1 \cup
%\boldsymbol{\beta}_1$ whose belt circle is the given curve.

Let $P'= S \cap f^{-1}(3/2).$ Then $\partial P'$ is the union of
$\partial P$ and some of the feet of the additional tubes. Next we
are going to modify $P'$ such that it becomes disjoint from these
additional tubes and it defines a surface equivalent to $S$.

Let $S_0$ be a component of $S$ and let $P_0' = P' \cap S_0.$ Since
$\partial S_0$ intersects both $R_+(\gamma)$ and $R_-(\gamma)$ we
see that $A \cap P_0' \neq \emptyset$ and $B \cap P_0' \neq
\emptyset.$ Because $S_0$ is connected $P_0'$ is also connected.
Note that for $\alpha \in \boldsymbol{\alpha}_1$ we have $\alpha
\cap P' = \emptyset.$ Thus we can achieve using isotopies that every
arc of $\alpha \cap P'$ for each $\alpha \in \boldsymbol{\alpha}_0$
intersects $A.$ Indeed, for every component $P_0'$ of $P'$ choose an
arc $\varphi_0 \subset P_0'$ whose endpoint lies on $A$ and
intersects every $\alpha$-arc lying in $P'_0.$ Then simultaneously
apply a finger move along $\varphi_0$ to all the $\alpha$-arcs that
intersect $\varphi_0.$ Similarly, we can achieve that each arc of
$\beta \cap P'$ for every $\beta \in \boldsymbol{\beta}_0$
intersects $B.$ This can be done keeping both the $\alpha$- and the
$\beta$-curves pairwise disjoint.
%This can be done without introducing any
%intersection points between both the $\alpha$- and the
%$\beta$-curves.

Let $F \subset \partial P'$ be the foot of a tube whose belt circle
is a curve $\alpha_1 \in \boldsymbol{\alpha}_1.$ Pick a point $p \in
F.$ Since every arc of $\beta \cap P'$ for $\beta \in
\boldsymbol{\beta}_0$ intersects $B$ each component of $P' \setminus
(\cup \boldsymbol{\beta}_0)$ intersects $B.$ Thus we can connect $p$
to $B$ with an arc $\varphi$ lying in $P' \setminus (\cup
\boldsymbol{\beta}).$ Now handleslide every $\alpha \in
\boldsymbol{\alpha}_0$ that intersects $\varphi$ over $\alpha_1$
along $\varphi.$ Then we can handleslide $B$ over $\alpha_1$ along
$\varphi.$ To this handleslide corresponds an isotopy of $S$ through
decomposing surfaces such that $S \cap f^{-1}(3/2)$ changes the
required way (given by taking the negative gradient flow lines of
$f$ flowing out of $B$). Thus we have removed $F$ from $P'.$ The
case when the belt circle of the tube lies in $\boldsymbol{\beta}_1$
is completely analogous. By repeating this process we can remove all
the additional one-handles from $P'.$ Call this new quasi-polygon
$P.$

%Finally, delete every $\alpha$ curve such that the corresponding
%index one critical point cancels with an index zero critical point
%and delete every $\beta$ curve such that the corresponding index two
%critical point cancels with an index three critical point.

Finally, cancel every index zero critical point with an index one
critical point and every index three critical point with an index
two critical point and delete the corresponding $\alpha$- and
$\beta$-curves. The balanced diagram obtained this way, together
with the quasi-polygon $P,$ defines $S.$
\end{proof}

\begin{lem} \label{lem:2}
Let $(M,\gamma) \rightsquigarrow^S (M',\gamma')$ be a surface
decomposition such that for every component $V$ of $R(\gamma)$ the
set of closed components of $S \cap V$ consists of parallel oriented
boundary-coherent simple closed curves. Then $S$ is isotopic to a
decomposing surface $S'$ such that each component of $\partial S'$
intersects both $R_+(\gamma)$ and $R_-(\gamma)$ and decomposing
$(M,\gamma)$ along $S'$ also gives $(M',\gamma').$ Furthermore, $O_S
= O_{S'}.$
\end{lem}

\begin{proof}
We call a tangency between two curves positive if their positive
unit tangent vectors coincide at the tangency point. Our main
observation is the following. Isotope a small arc of $\partial S$ on
$\partial M$ using a finger move through $\gamma$ such that during
the isotopy we have a positive tangency between $\partial S$ and
$s(\gamma)$ (thus introducing two new intersection points between
$\partial S$ and $s(\gamma)$). Let the resulting isotopy of
$\partial S$ be $\{\,s_t \colon 0\le t \le 1\,\}.$ Attach the collar
$\partial M \times I$ to $M$ to get a new manifold $\widetilde{M}$
and attach $\cup_{t \in I} (s_t \times \{t\})$ to $S$ to obtain a
surface $\widetilde{S} \subset \widetilde{M}.$ Then decomposing
$$(\widetilde{M}, \gamma \times \{1\})\approx (M, \gamma)$$ along
$\widetilde{S}$ we also get $(M',\gamma'),$ see Figure \ref{fig:5}.
Furthermore, $\widetilde{S}$ is isotopic to $S.$

Let $\gamma_0$ be a component of $\gamma$ such that $\gamma_0 \cap
\partial S$ consists of closed curves $\sigma_1,\dots, \sigma_k$. First
isotope $S$ in a neighborhood of $\partial S \cap \gamma_0$ through
decomposing surfaces such that after the isotopy $\sigma_1, \dots,
\sigma_k$ are all parallel to $s(\gamma)$ and $\nu_S$ points out of
$M$ along $\partial S \cap \gamma_0.$ This new decomposing surface
is equivalent to the original. Then isotope $\sigma_1, \dots,
\sigma_k$ into $R_-(\gamma).$ Decomposing along $S$ still gives
$(M',\gamma').$ Let $\delta$ be an oriented arc that intersects
$\sigma_1 \dots, \sigma_k,$ and $s(\gamma)$ exactly once and its
endpoint lies in $R_+(\gamma).$ Applying a finger move to $\sigma_1,
\dots, \sigma_k$ simultaneously along $\delta$ we get a positive
tangency between each $\sigma_i$ and $s(\gamma)$ since they are
oriented coherently.

%Choose two distinct points $p,q \in s(\gamma_0)$ and an arc $J_p
%\subset \gamma_0$ connecting $p$ and $R_+(\gamma)$ and an arc $J_q
%\subset \gamma_0$ connecting $q$ and $R_-(\gamma).$ Apply a finger
%move to $\gamma_0 \cap S$ along $J_p \cup J_q.$ Then every component
%of $\gamma_0 \cap S$ will intersect $s(\gamma).$ During the isotopy
%we only get positive tangencies because components of $\gamma_0 \cap
%S$ are oriented coherently with $s(\gamma_0).$ Now shrink $\gamma_0$
%such that $\gamma_0 \cap S$ only consists of essential arcs
%intersecting $s(\gamma)$ once.

Let $V$ be a component of $R(\gamma)$ and let $C_1, \dots, C_k$ be
the parallel oriented closed components of $S \cap V.$ Choose a
small arc $T$ that intersects every $C_i$ in a single point. Let
$\partial T = \{x,y\}.$ First suppose that $[C_1] \neq 0$ in
$H_1(V;\mathbb{Z}).$ Then we can connect both $x$ and $y$ to
$s(\gamma)$ by an arc whose interior lies in $\partial M \setminus
(\partial S \cup s(\gamma)).$ This is possible since $C_1$ does not
separate $\partial V$ and now $\partial S \cap \gamma$ has no closed
components. This way we obtain an arc $\delta \subset
\partial M$ such that for every $1 \le i \le k$ we have $|\delta
\cap C_i| = 1$ and $\partial \delta = \delta \cap s(\gamma);$
moreover, $$\delta \cap \partial S = \delta \cap (C_1 \cup \dots
\cup C_k).$$ Recall that $s(\gamma)$ is oriented coherently with
$\partial V$ (this is especially important if $s(\gamma)$ is
disconnected and $\delta$ connects two distinct components of
$s(\gamma)$) and the curves $C_1, \dots, C_k$ are also oriented
coherently. Thus with exactly one of the orientations of $\delta$ if
we apply a finger move to all the $C_i$ simultaneously we get a
positive tangency between each $C_i$ and $\partial V,$ and thus also
$s(\gamma).$

Now suppose that $[C_1] = 0$ in $H_1(V;\mathbb{Z})$ and $C_1$ is
oriented as the boundary of its interior. Then exactly one of $x$
and $y$ can be connected to $s(\gamma)$ by an arc $\delta_0$ whose
interior lies in $\partial M \setminus (\partial S \cup s(\gamma)).$
The arc $T \cup \delta_0$ defines an oriented arc $\delta$ whose
endpoint lies on $s(\gamma).$ If we apply a finger move to each
$C_i$ along $\delta$ we get positive tangencies with $s(\gamma)$
because every $C_i$ is oriented as the boundary of its interior and
$s(\gamma)$ is oriented coherently with respect to $\partial V$.

Continuing this process we get a surface $S'$ isotopic to $S$ such
that each component of $\partial S'$ intersects $s(\gamma)$ and
decomposing $(M, \gamma)$ along $S'$ we still get $(M', \gamma').$

\begin{figure}[t]
\includegraphics{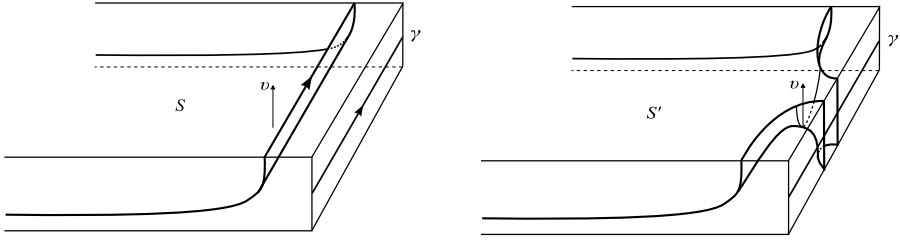}
\caption{Making a decomposing surface good.} \label{fig:5}
\end{figure}

To show that $O_S = O_{S'}$ first observe that if $S_0$ and $S_1$
are equivalent then $O_{S_0} = O_{S_1}.$ Now suppose that for some
component $\gamma_0$ of $\gamma$ the components of $\partial S \cap
\gamma_0$ are curves $\sigma_1, \dots, \sigma_k$ parallel to
$s(\gamma)$ such that $\nu_S$ points out of $M$ along them.
Moreover, suppose that $S'$ only differs from $S$ by isotoping
$\sigma_1,\dots, \sigma_k$ into $R_-(\gamma).$ If $\mathfrak{s}$ is
a $\text{Spin}^c$ structure and $v$ is a vector field representing
it, then in a standard neighborhood of $\gamma_0$ we have $v \neq
\pm\nu_S$ and $v \neq \pm\nu_{S'}.$ So $\mathfrak{s} \in O_S$ if and
only if $\mathfrak{s} \in O_{S'}.$

Thus we only have to show that $O_S = O_{S'}$ when $S$ and $S'$ are
related by a small finger move of $\partial S$ that crosses
$s(\gamma)$ through a positive tangency. Let $\mathfrak{s}$ be a
$\text{Spin}^c$ structure on $(M, \gamma)$ and $v$ a vector field
representing it. Then in a standard neighborhood $U$ of the tangency
point we can perform the isotopy such that in $U$ we have $v \neq
\pm\nu_S;$ furthermore, $v^{\perp}$ and $S'$ only have a single
hyperbolic tangency, where $v = \nu_S$ (see Figure \ref{fig:5}).
Thus $\mathfrak{s} \in O_S$ if and only if $\mathfrak{s} \in
O_{S'}.$ Note that if the tangency of $\partial S$ and $s(\gamma)$
is negative during the isotopy then at the hyperbolic tangency $v =
-\nu_{S'}.$

If $(M,\gamma)$ is strongly balanced then $O_S = O_{S'}$ also
follows from Lemma \ref{lem:1}. Indeed, $\left\langle\,
c_1(\mathfrak{s},t), [S] \,\right\rangle$ is invariant under
isotopies of $S.$ As before, we can suppose that the closed
components of $\partial S \cap \gamma$ are parallel to $s(\gamma)$
and $\nu_S$ points out of $M$ along them. In the above proof $I$ and
$r$ are unchanged when we isotope $\sigma_i$ from $\gamma_0$ to
$R_-(\gamma)$ since we can achieve that $\nu_S$ and $v$ are never
parallel along $\partial S,$ so $I$ and $r$ change continuously.
When we do a finger move $I$ decreases by $1$ according to part (1)
of Lemma \ref{lem:9} and $r$ also decreases by $1,$ as can be seen
from Figure \ref{fig:5}. Thus $c(S,t)=c(S',t).$
\end{proof}

\begin{defn} \label{defn:17}
We call a decomposing surface $S \subset (M,\gamma)$ \emph{good} if
it is open and each component of $\partial S$ intersects both
$R_+(\gamma)$ and $R_-(\gamma).$ We call a surface diagram $(\Sigma,
\boldsymbol{\alpha}, \boldsymbol{\beta},P)$ \emph{good} if $A$ and
$B$ have no closed components.
\end{defn}

\begin{rem} \label{rem:3}
Because of Lemma \ref{lem:2} it is sufficient to prove Theorem
\ref{thm:1} for good decomposing surfaces. According to Proposition
\ref{prop:4} for each good decomposing surface we can find a good
surface diagram adapted to it.
\end{rem}

\begin{prop} \label{prop:9}
Suppose that $S$ is a good decomposing surface in the balanced
sutured manifold $(M,\gamma).$ Then there exists an
\emph{admissible} surface diagram of $(M,\gamma)$ adapted to $S.$
\end{prop}

\begin{proof}
According to Remark \ref{rem:3} we can find a good surface diagram
$(\Sigma, \boldsymbol{\alpha}, \boldsymbol{\beta}, P)$ adapted to
$S.$

Here we improve on the idea of the proof of \cite[Proposition
3.15]{sutured}. Choose pairwise disjoint arcs $\gamma_1, \dots,
\gamma_k \subset \Sigma \setminus B$ whose endpoints lie on
$\partial \Sigma$ and together generate $H_1(\Sigma \setminus
B,\partial(\Sigma \setminus B); \mathbb{Z}).$ This is possible
because each component of $\partial (\Sigma \setminus B)$ intersects
$\partial \Sigma.$ Choose curves $\gamma_1', \dots, \gamma_k'$ such
that $\gamma_i$ and $\gamma_i'$ are parallel and oriented
oppositely.

Then wind the $\boldsymbol{\alpha}$ curves along $\gamma_1,
\gamma_1', \dots, \gamma_k, \gamma_k'$ as in the proof of
\cite[Proposition 3.15]{sutured}. A similar argument as there gives
that after the winding $(\Sigma, \boldsymbol{\alpha},
\boldsymbol{\beta})$ will be admissible. Note that every $\alpha \in
\boldsymbol{\alpha}$ lies in $\Sigma \setminus B.$ Thus if a linear
combination $\mathcal{A}$ of $\alpha$-curves intersects every
$\gamma_i$ algebraically zero times then $\mathcal{A}$ is
null-homologous in $\Sigma \setminus B,$ and thus also in $\Sigma.$
Since the winding is done away from $B$ the new diagram is still
adapted to $S.$
\end{proof}

\section{Balanced diagrams and surface decompositions}

\begin{defn} \label{defn:18}
Let $(\Sigma, \boldsymbol{\alpha}, \boldsymbol{\beta}, P)$ be a
surface diagram (see Definition \ref{defn:16}). Then we can uniquely
associate to it a tuple $D(P) = (\Sigma', \boldsymbol{\alpha}',
\boldsymbol{\beta}', P_A, P_B, p),$ where $(\Sigma',
\boldsymbol{\alpha}', \boldsymbol{\beta}')$ is a balanced diagram,
$p \colon \Sigma' \to \Sigma$ is a smooth map, and $P_A, P_B \subset
\Sigma'$ are two closed subsurfaces (see Figure \ref{fig:2}).

To define $\Sigma'$ take two disjoint copies of $P$ that we call
$P_A$ and $P_B$ together with diffeomorphisms $p_A \colon P_A \to P$
and $p_B \colon P_B \to P.$ Cut $\Sigma$ along $\partial P$ and
remove $P.$ Then glue $A$ to $P_A$ using $p_A^{-1}$ and $B$ to $P_B$
using $p_B^{-1}$ to obtain $\Sigma'.$ The map $p \colon \Sigma' \to
\Sigma$ agrees with $p_A$ on $P_A$ and $p_B$ on $P_B,$ and it maps
$\Sigma' \setminus (P_A \cup P_B)$ to $\Sigma \setminus P$ using the
obvious diffeomorphism. Finally, let $\boldsymbol{\alpha}' =
\{\,p^{-1}(\alpha) \setminus P_B \colon \alpha \in
\boldsymbol{\alpha}\,\}$ and $\boldsymbol{\beta}' =
\{\,p^{-1}(\beta) \setminus P_A \colon \beta \in
\boldsymbol{\beta}\,\}.$

%We obtain $\Sigma'$  as follows. Let $\Sigma_1$ be the disjoint
%union of $\Sigma$ and a copy $P_B$ of $P$ and let $p_B \colon P_B
%\to P$ be a homeomorphism. Cut $\Sigma_1$ along $B$ to obtain a
%surface $\Sigma_2.$ Denote by $b_1, \dots, b_m$ the components of
%$B.$ For every $1 \le i \le m$ let $b_i'$ be the image of $b_i$ in
%$\Sigma \setminus \text{int}(P).$ After gluing $b_i'$ to
%$p_B^{-1}(b_i)$ for every $1 \le i \le m$ we obtain the surface
%$\Sigma'.$ The image of the quasi-polygon $P$ in $\Sigma'$ is called
%$P_A.$ Then $p \colon \Sigma' \to \Sigma$ is defined by $p|P_B =
%p_B$ and $p|(\Sigma' \setminus P_B) = \text{Id}.$ Finally, let
%$\boldsymbol{\alpha}' = \{\,p^{-1}(\alpha) \setminus P_B \colon
%\alpha \in \boldsymbol{\alpha}\,\}$ and $\boldsymbol{\beta}' =
%\{\,p^{-1}(\beta) \setminus P_A \colon \beta \in
%\boldsymbol{\beta}\,\}.$

$D(P)$ is uniquely characterized by the following properties. The
map $p$ is a local diffeomorphism in $\text{int}(\Sigma');$
furthermore, $p^{-1}(P)$ is the disjoint union of $P_A$ and $P_B.$
Moreover, $p|P_A \colon P_A \to P,$ and $p|P_B \colon P_B \to P,$
and also
$$p|(\Sigma' \setminus (P_A \cup P_B)) \colon \Sigma' \setminus (P_A
\cup P_B) \to \Sigma \setminus P$$ are diffeomorphisms. Furthermore,
$p(\text{int}(\Sigma') \cap \partial P_A) = \text{int}(A)$ and
$p(\text{int}(\Sigma') \cap \partial P_B) = \text{int}(B).$ Finally,
$p|(\cup\boldsymbol{\alpha}') \colon \cup\boldsymbol{\alpha}' \to
\cup\boldsymbol{\alpha}$ and $p|(\cup\boldsymbol{\beta}') \colon
\cup\boldsymbol{\beta}' \to \cup\boldsymbol{\beta}$ are
diffeomorphisms. Thus $(\cup\boldsymbol{\alpha}') \cap P_B =
\emptyset$ and $(\cup\boldsymbol{\beta}') \cap P_A = \emptyset.$

\begin{figure}[t]
\includegraphics{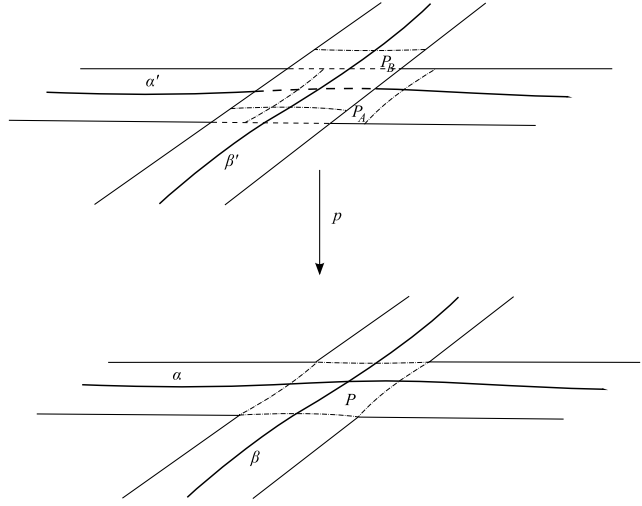}
\caption{Balanced diagrams before and after a surface
decomposition.} \label{fig:2}
\end{figure}

There is a unique holomorphic structure on $\Sigma'$ that makes the
map $p$ holomorphic. Since $p$ is a local diffeomorphism in
$\text{int}(\Sigma)$ it is even conformal.
\end{defn}

So $p$ is $1:1$ over $\Sigma \setminus P,$ it is $2:1$ over $P,$ and
$\alpha$ curves are lifted to $P_A$ and $\beta$ curves to $P_B.$

\begin{prop} \label{prop:5}
Let $(M,\gamma)$ be a balanced sutured manifold and $$(M,
\gamma)\rightsquigarrow^S (M',\gamma')$$ a surface decomposition. If
$(\Sigma, \boldsymbol{\alpha},\boldsymbol{\beta},P)$ is a surface
diagram adapted to $S$ and if $D(P) = (\Sigma',
\boldsymbol{\alpha}', \boldsymbol{\beta}', P_A, P_B, p)$ then
$(\Sigma', \boldsymbol{\alpha}', \boldsymbol{\beta}')$ is a balanced
diagram defining $(M',\gamma').$
\end{prop}

\begin{proof}
Let $(M_1,\gamma_1)$ be the sutured manifold defined by the diagram
$(\Sigma', \boldsymbol{\alpha}', \boldsymbol{\beta}').$ We are going
to construct an orientation preserving homeomorphism $h \colon
(M_1,\gamma_1) \to (M', \gamma')$ that takes $R_+(\gamma_1)$ to
$R_+(\gamma').$ Figure \ref{fig:4} is a schematic illustration of
the proof.

Let $N_A$ and $N_B$ be regular neighborhoods of $P_A$ and $P_B$ in
$\Sigma'$ so small that $\alpha' \cap N_B = \emptyset$ and $\beta'
\cap N_A = \emptyset$ for every $\alpha' \in \boldsymbol{\alpha}'$
and $\beta' \in \boldsymbol{\beta}'.$ Furthermore, let $N = N_A \cup
N_B.$ Define $\lambda \colon \Sigma' \to I$ to be a smooth function
such that $\lambda(x) = 1$ for $ x \in \Sigma' \setminus N$ and
$\lambda(x) =1/2$ for $x \in P_A \cup P_B.$ Moreover, let $\mu
\colon \Sigma' \to I$ be a smooth function such that $\mu(x) = 1 -
\lambda(x)$ for $x \in N_B$ and $\mu(x) = 0$ for $x \in \Sigma'
\setminus N_B.$

The homeomorphism $h$ is constructed as follows. For $(x,t) \in
\Sigma' \times I$ let $$h(x,t) = \left( p(x), \mu(x) + \lambda(x)t
\right).$$ Since for every $x \in \Sigma'$ and $t \in I$ the
inequality $0 \le \mu(x) + \lambda(x)t \le 1$ holds the map $h$
takes $\Sigma' \times I$ into $\Sigma \times I \subset (M,\gamma).$
Choose an $\alpha' \in \boldsymbol{\alpha}'$ and let $\alpha =
p(\alpha') \in \boldsymbol{\alpha}.$ Let $D_{\alpha'}$ be the
2-handle attached to $\Sigma' \times I$ along $\alpha' \times \{0\}$
and $D_{\alpha}$ the 2-handle attached to $\Sigma \times I$ along
$\alpha \times \{0\}.$ Since $\alpha' \cap N_B = \emptyset$ and
because $\mu(x) + \lambda(x) \cdot 0 =0$ for $x \in \Sigma'
\setminus N_B$ we see that $h(\alpha' \times \{0\}) = \alpha \times
\{0\}.$ Thus $h$ naturally extends to a map from $(\Sigma' \times I)
\cup D_{\alpha'}$ to $(\Sigma \times I) \cup D_{\alpha}.$ Similarly,
for $\beta' \in \boldsymbol{\beta}'$ we have $\beta' \cap N_A =
\emptyset.$ Furthermore, $\mu(x) + \lambda(x) \cdot 1 = 1$ for $x
\in \Sigma' \setminus N_A.$ Thus $h$ also extends to the 2-handles
attached along the $\beta$-curves. So now we have a local
homeomorphism from $(M_1,\gamma_1)$ into $(M,\gamma).$

Recall that $S \subset (M, \gamma)$ is equivalent to the surface
obtained by smoothing $$(P \times \{1/2\}) \cup (A \times [1/2,1])
\cup (B \times [0,1/2]) \subset \Sigma \times I.$$ Since $h(N_B
\times \{0\}) \cup h(N_A \times \{1\})$ is a smoothing of the above
surface we can assume that it is in fact equal to $S.$ Indeed, for
$x \in P_A$ we have that $\mu(x) + \lambda(x) \cdot 1 = 1/2$ and for
$x \in P_B$ the equality $\mu(x) + \lambda(x) \cdot 0 = 1-\lambda(x)
= 1/2$ holds. Moreover, $p(\partial N_A \setminus
\partial \Sigma') = A'$ is a curve parallel to $A,$ thus for $x \in
\partial N_A \setminus \partial \Sigma'$ we have $h(x,1) \in A' \times
\{1\}.$ Similarly, $h(x,0) \in B' \times \{0\}$ for $x \in \partial
N_B \setminus \partial \Sigma',$ where $B'$ is a curve parallel and
close to $B.$

\begin{figure}[t]
\includegraphics{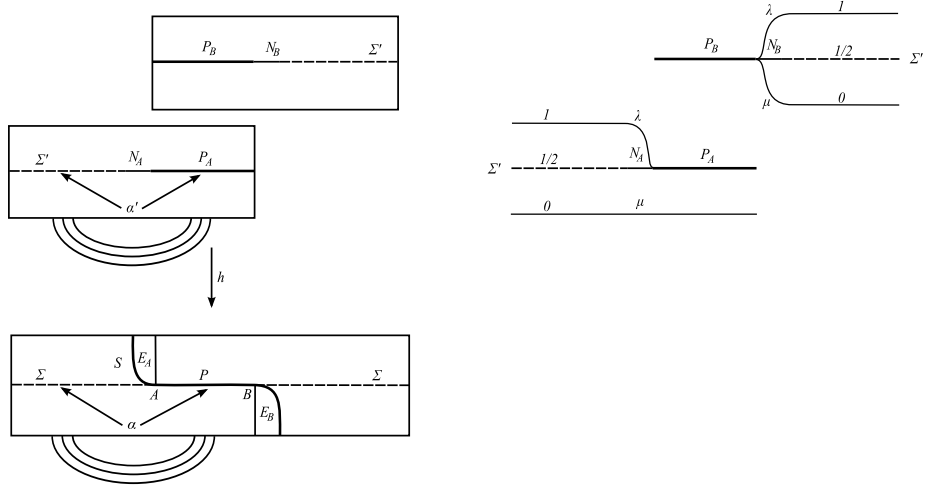}
\caption{The left hand side shows the homeomorphism $h.$ On the
right we can see the functions $\lambda$ and $\mu.$} \label{fig:4}
\end{figure}

Let $E_A \subset \Sigma \times I$ be the set of points $(y,s)$ such
that $y = p(x)$ for some $x \in N_A \setminus P_A$ and $s \ge \mu(x)
+ \lambda(x).$ Define $E_B \subset \Sigma \times I$ to be the set of
those points $(y,s)$ such that $y = p(x)$ for some $x \in N_B
\setminus P_B$ and $s \le \mu(x).$ Now we are going to show that the
map $$h | (\Sigma' \times I \setminus (P_A \times \{1\} \cup P_B
\times \{0\})) \to (\Sigma \times I) \setminus (S \cup E_A \cup
E_B)$$ is a homeomorphism by constructing its continuous inverse.
Let $$(y,s) \in (\Sigma \times I) \setminus (S \cup E_A \cup E_B).$$
If $y \in \Sigma \setminus p(N)$ then $h^{-1}(y,s) = (p^{-1}(y),s).$
If $y \in P$ and $s < 1/2$ then $h^{-1}(y,s) = (p^{-1}(y) \cap P_A,
2s)$ and for $s > 1/2$ we have $h^{-1}(y,s) = (p^{-1}(y) \cap P_B,
2s-1).$ In the case when $y \in p(N_A \setminus P_A)$ and $s <
\mu(x)+\lambda(x)$ we let $h^{-1}(y,s) = (x,t),$ where $x =
p^{-1}(y)$ and $t = (s-\mu(x))/\lambda(x) < 1.$ Note that here
$\mu(x) = 0,$ and thus $t \ge 0.$ Finally, for $y \in p(N_B
\setminus P_B)$ and $s > \mu(x)$ define $h(y,s) = (x,t),$ where $x =
p^{-1}(y)$ and $t = (s -\mu(x))/\lambda(x) > 0.$ Here $t \le 1$
because $s \le 1$ and $\mu(x) = 1 - \lambda(x).$

Recall that we defined the surfaces $S'_+$ and $S'_-$ in Definition
\ref{defn:4}. Since $S$ is oriented coherently with $P \times
\{1/2\}$ thickening $S'_+ \cap R_-(\gamma)$ in $\partial M'$ can be
achieved by cutting off its neighborhood $E_B$ and taking $B \times
[0,1/2] \subset \partial E_B$ to belong to $\gamma'.$ Similarly,
$E_A$ is a neighborhood of $S'_- \cap R_+(\gamma)$ in $M',$ and
cutting it off from $M'$ we can add $A \times [1/2,1]$ to $\gamma'.$
Thus we can identify $M'$ with the metric completion of $M \setminus
(S \cup E_A \cup E_B)$ and $\gamma'$ with $(\gamma \cap M') \cup (A
\times [1/2,1]) \cup (B \times [0,1/2]).$

What remains is to show that $h(\gamma_1) = \gamma'.$ If $x \in
(\partial \Sigma') \setminus (P_A \cup P_B)$ then for any $t \in I$
we have $$h(x,t) = (p(x), \mu(x) + \lambda(x)t) \in \gamma \cap M'
\subset \gamma'$$ because $p(x) \in \partial \Sigma.$ On the other
hand, for $x \in \partial \Sigma' \cap P_A$ and $t \in I$ we have
$h(x,t) \in B \times [0,1/2],$ which is part of $\gamma'$ by the
above construction. The case $x \in \partial \Sigma' \cap P_B$ is
similar.
\end{proof}

\begin{defn} \label{defn:19}
Let $(\Sigma, \boldsymbol{\alpha}, \boldsymbol{\beta}, P)$ be a
surface diagram. We call an intersection point $\mathbf{x} \in
\mathbb{T}_{\alpha} \cap \mathbb{T}_{\beta}$ \emph{outer} if
$\mathbf{x} \cap P = \emptyset.$ We denote by $O_P$ the set of outer
intersection points. Then $I_P = (\mathbb{T}_{\alpha} \cap
\mathbb{T}_{\beta}) \setminus O_P$ is called the set of \emph{inner}
intersection points.
\end{defn}

\begin{lem} \label{lem:3}
Let $(M, \gamma) \rightsquigarrow^S (M', \gamma')$ be a surface
decomposition and suppose that $(\Sigma, \boldsymbol{\alpha},
\boldsymbol{\beta}, P)$ is a surface diagram adapted to $S.$ Let
$\mathbf{x} \in \mathbb{T}_{\alpha} \cap \mathbb{T}_{\beta}.$ Then
$\mathbf{x} \in O_P$ if and only if $\mathfrak{s}(\mathbf{x}) \in
O_S.$ Furthermore, if $D(P) = (\Sigma', \boldsymbol{\alpha}',
\boldsymbol{\beta}', P_A, P_B, p)$ then $p$ gives a bijection
between $\mathbb{T}_{\alpha'} \cap \mathbb{T}_{\beta'}$ and $O_P.$
\end{lem}

\begin{proof}
Let $f$ be a Morse function on $M$ compatible with the diagram
$(\Sigma, \boldsymbol{\alpha}, \boldsymbol{\beta}).$ If $\mathbf{x}
\in O_P$ then the multi-trajectory $\gamma_{\mathbf{x}}$ (see
Definition \ref{defn:20}) is disjoint from $S.$ Consequently, the
regular neighborhood $N(\gamma_{\mathbf{x}})$ can be chosen to be
disjoint from $S.$ Thus $\mathfrak{s}(\mathbf{x})$ can be
represented by a unit vector field $v$ that agrees with
$\text{grad}(f)/\norm{\text{grad}(f)}$ in a neighborhood of $S.$
Since the orientation of $S$ is compatible with the orientation of
$P \subset \Sigma,$ even after smoothing the corners of $(P\times
\{1/2\}) \cup (A \times [1/2,1]) \cup (B \times [0,1/2])$ we have
that $v$ is nowhere equal to $-\nu_S.$ So we see that
$\mathfrak{s}(\mathbf{x}) \in O_S.$

\begin{figure}[t]
\includegraphics{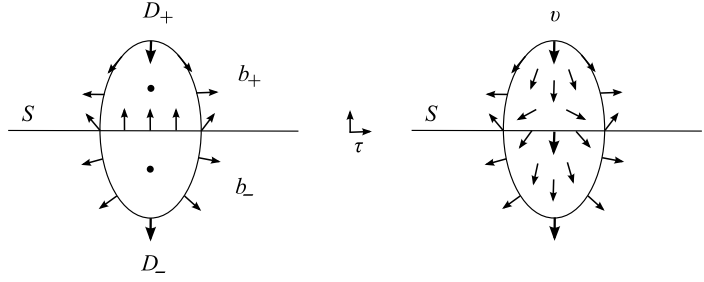}
\caption{This is a schematic two-dimensional picture illustrating
the proof of Lemma \ref{lem:3}.} \label{fig:8}
\end{figure}

Now suppose that $\mathbf{x} \in I_P.$ Let $\gamma_{\mathbf{x}}$ be
the multi-trajectory associated to $\mathbf{x}.$ Since $S$ is open
its tangent bundle $TS$ is trivial. Thus there is a trivialization
$\tau = (\tau_1, \tau_2, \tau_3)$ of $TM|(S \cup
N(\gamma_{\mathbf{x}}))$ such that $\tau_3|S = \nu_S$ and
$(\tau_1|S, \tau_2|S)$ is a trivialization of $TS.$ The
$\text{Spin}^c$ structure $\mathfrak{s}(\mathbf{x})$ can be
represented by a unit vector field $v$ such that $v|(M \setminus
N(\gamma_{\mathbf{x}}))$ agrees with
$$g = \frac{\text{grad}(f)|(M \setminus
N(\gamma_{\mathbf{x}}))}{\norm{\text{grad}(f)|(M \setminus
N(\gamma_{\mathbf{x}}))}}.$$ If $v$ was outer then for any ball $B^3
\subset M \setminus S$ the vector field $v|(M \setminus B^3)$ would
be homotopic through unit vector fields rel $\partial M$ to a field
$v'$ such that $v'|S$ is nowhere equal to $-\nu_S.$ So to prove that
$\mathfrak{s}(\mathbf{x}) \not\in O_S$ it is sufficient to show that
$v|S$ is not homotopic through unit vector fields rel $\partial S$
to a vector field $v'$ on $S$ that is nowhere equal to $-\nu_S.$ In
the trivialization $\tau$ we can think of $v|(S \cup
N(\gamma_{\mathbf{x}}))$ as a map from $S \cup
N(\gamma_{\mathbf{x}})$ to $S^2$ and $-\nu_S$ corresponds to the
South Pole $s \in S^2.$ If we put $S$ in generic position $v_0 =
v|\partial M$ is nowhere equal to $-\nu_S.$
%(we can even suppose that $v_0| \partial S$ is tangent to $S$).
Thus $v$ maps $\partial S$ into $S^2 \setminus \{s\}.$

Let $x \in \mathbf{x}$ and let $\gamma_x$ be the component of
$\gamma_{\mathbf{x}}$ containing $x.$ Then $\gamma_x \cap S =
\emptyset$ if $x \not \in P$ and $\gamma_x \cap S = \{x\}$ if $x \in
P.$ So suppose that $x \in P.$ We denote $N(\gamma_x)$ by $B$ and
let $B_+$ and $B_-$ be the closures of the two components of $B
\setminus S;$ an index one critical point of $f$ lies in $B_-$ and
an index two critical point in $B_+.$ Moreover, let $D_{\pm} =
\partial B_{\pm} \setminus S.$ The vector field
$\text{grad}(f)|B$ is a map from $B$ to $\mathbb{R}^3$ in the
trivialization $\tau.$ Let $$b_{\pm} = \frac{\text{grad}(f)|\partial
B_{\pm}}{\norm{\text{grad}(f)|\partial B_{\pm}}},$$ see Figure
\ref{fig:8}. Since $B_{\pm}$ contains an index $\pm 1$ singularity
of $\text{grad}(f)$ we see that $\#b_{\pm}^{-1}(s) = \pm 1.$ Here
$\#$ denotes the algebraic number of points in a given set. Since
$\text{grad}(f)|(S \cap B)$ is equal to $\nu_S$ we even get that
$\#(b_{\pm}^{-1}(s) \cap D_{\pm}) = \pm 1.$ Let $v_{\pm} =
v|\partial B_{\pm}.$ Then $\#v_{\pm}^{-1}(s) = 0$ because $v$ is
nowhere zero. The co-orientation of $S$ is given by
$\text{grad}(f),$ so $S \cap B \subset S$ is oriented coherently
with $\partial B_-.$ Moreover, $v| D_- = b_-|D_-,$ so we see that
$\#(v_-^{-1}(s) \cap S) = 1.$ We have seen that $g|(S \setminus P) =
v|(S \setminus P)$ is nowhere equal to $-\nu_S.$ So we conclude that
$\#(v|S)^{-1}(s) = |\mathbf{x} \cap P|.$ Thus if $\mathbf{x} \in
I_P$ then $v|S$ is not homotopic to a map $S \to S^2 \setminus
\{s\}$ through a homotopy fixing $\partial S.$ This means that
$\mathfrak{s}(x) \not\in O_S.$

The last part of the statement follows from the fact that $p$ is a
diffeomorphism between $\Sigma' \setminus (P_A \cup P_B)$ and
$\Sigma \setminus P,$ furthermore $(\cup\boldsymbol{\alpha}') \cap
P_B = \emptyset$ and $(\cup\boldsymbol{\beta}') \cap P_A =
\emptyset.$
\end{proof}

\begin{rem} \label{rem:5}
We can slightly simplify the proof of Lemma \ref{lem:3} when $O_P
\neq \emptyset.$ Suppose that $\mathbf{x} \in I_P$ and let
$\mathbf{y} \in O_P$ be an arbitrary intersection point. Using
\cite[Lemma 4.7]{sutured} we get that $\mathfrak{s}(\mathbf{x}) -
\mathfrak{s}(\mathbf{y}) = PD[\gamma_{\mathbf{x}} -
\gamma_{\mathbf{y}}].$ Since the co-orientation of $P \subset S$ is
given by $\text{grad}(f)$ we get that $$\langle
\,\mathfrak{s}(\mathbf{x}) - \mathfrak{s}(\mathbf{y}), [S] \,\rangle
= |\gamma_{\mathbf{x}} \cap S| - |\gamma_{\mathbf{y}} \cap S| =
|\mathbf{x} \cap P| - |\mathbf{y} \cap P| \neq 0.$$ If
$\mathfrak{s}(\mathbf{x})$ was outer then both
$\mathfrak{s}(\mathbf{x})$ and $\mathfrak{s}(\mathbf{y})$ could be
represented by unit vector fields that are homotopic over $S$ rel
$\partial S$ since $(STM|S) \setminus (-\nu_S)$ is a bundle with
contractible fibers. And that would imply that $\langle
\,\mathfrak{s}(\mathbf{x}) - \mathfrak{s}(\mathbf{y}), [S] \,\rangle
= 0.$ Thus $\mathfrak{s}(x)$ is not outer.
\end{rem}

\begin{note}
We will also denote by $O_P$ and $I_P$ the subgroups of $CF(\Sigma,
\boldsymbol{\alpha}, \boldsymbol{\beta})$ generated by the outer and
inner intersection points, respectively.
\end{note}

\begin{cor} \label{cor:1}
For a surface diagram $(\Sigma, \boldsymbol{\alpha},
\boldsymbol{\beta}, P)$ such that $(\Sigma, \boldsymbol{\alpha},
\boldsymbol{\beta})$ is admissible the chain complex $(CF(\Sigma,
\boldsymbol{\alpha}, \boldsymbol{\beta}),
\partial)$ is the direct sum of the subcomplexes $(O_P, \partial |
O_P)$ and $(I_P, \partial | I_P).$
\end{cor}

\section{An algorithm providing a nice surface diagram}

In this section we generalize the results of \cite{Sucharit} to
sutured Floer homology and surface diagrams. Our argument is an
elaboration of the Sarkar-Wang algorithm. The basic approach is the
same, but there are some important differences. The definition of
distance had to be modified to work in this generality. Additional
technical difficulties arise because when we would like to make a
surface diagram nice we have to assure that the property $A \cap B =
\emptyset$ is preserved. Moreover, $\boldsymbol{\alpha}$ or
$\boldsymbol{\beta}$ might not span $H_1(\Sigma;\mathbb{Z}),$ which
makes some of the arguments more involved.

\begin{defn} \label{defn:25}
We say that the surface diagram $(\Sigma, \boldsymbol{\alpha},
\boldsymbol{\beta}, P)$ is \emph{nice} if every component of $\Sigma
\setminus \left(\bigcup \boldsymbol{\alpha} \cup \bigcup
\boldsymbol{\beta} \cup A \cup B\right)$ whose closure is disjoint
from $\partial \Sigma$ is a bigon or a square. In particular, a
balanced diagram $(\Sigma, \boldsymbol{\alpha}, \boldsymbol{\beta})$
is called \emph{nice} if the surface diagram
$(\Sigma,\boldsymbol{\alpha}, \boldsymbol{\beta}, \emptyset)$ is
nice.
\end{defn}

\begin{defn} \label{defn:27}
Let $(\Sigma, \boldsymbol{\alpha}, \boldsymbol{\beta}, P)$ be a
surface diagram. Then a \emph{permissible move} is an isotopy or a
handle slide of the $\alpha$-curves in $\Sigma \setminus B$ or the
$\beta$-curves in $\Sigma \setminus A.$
\end{defn}

\begin{lem} \label{lem:4}
Let $\mathcal{S}$ be a surface diagram adapted to the decomposing
surface $S \subset (M, \gamma)$. If the surface diagram
$\mathcal{S}'$ is obtained from $\mathcal{S}$ using permissible
moves then $\mathcal{S}'$ is also adapted to $S.$
\end{lem}

\begin{proof}
This is a simple consequence of the definitions.
\end{proof}

\begin{thm} \label{thm:7}
Every good surface diagram $\mathcal{S}=(\Sigma,
\boldsymbol{\alpha}, \boldsymbol{\beta}, P)$ can be made nice using
permissible moves. If $(\Sigma, \boldsymbol{\alpha},
\boldsymbol{\beta})$ was admissible our algorithm gives an
admissible diagram.
\end{thm}

\begin{proof}
Let $\mathbb{A} = (\bigcup \boldsymbol{\alpha}) \cup B$ and
$\mathbb{B} = (\bigcup \boldsymbol{\beta}) \cup A.$ The set of those
components of $\Sigma \setminus (\mathbb{A} \cup \mathbb{B})$ whose
closure is disjoint from $\partial \Sigma$ is denoted by
$C(\mathcal{S}).$

First we achieve that every element of $C(\mathcal{S})$ is
homeomorphic to $D^2.$ Let $R(\mathcal{S})$ denote the set of those
elements of $C(\mathcal{S})$ which are \emph{not} homeomorphic to
$D^2$ and let $a(\mathcal{S}) = \sum_{R \in R(\mathcal{S})}
(1-\chi(R)).$ Choose a component $R \in R(\mathcal{S}).$ Then
$H_1(R,\partial R) \neq 0,$ thus there exists a curve
$(\delta,\partial \delta) \subset (R,
\partial R)$ such that $[\delta] \neq 0$ in $H_1(R,\partial R).$
Moreover, we can choose $\delta$ such that either $\delta(0) \in
\bigcup\boldsymbol{\alpha}$ and $\delta(1) \in \mathbb{B},$ or
$\delta(0) \in \bigcup{\boldsymbol{\beta}}$ and $\delta(1) \in
\mathbb{A},$ as follows. Since our surface diagram is good there are
no closed components of $A$ and $B,$ and note that $A \cap B =
\emptyset.$ Furthermore, $\partial R \cap \mathbb{A} \neq \emptyset$
and $\partial R \cap \mathbb{B} \neq \emptyset$ since otherwise $R$
would give a linear relation between either the $\alpha$-curves or
the $\beta$-curves. So if $\partial R$ is disconnected we can even
find two distinct components $C$ and $C'$ of $\partial R$ such that
$C \cap \mathbb{A} \neq \emptyset$ and $C' \cap \mathbb{B} \neq
\emptyset.$ Thus we can choose $\delta$ such that $\partial \delta
\cap \mathbb{A} \neq \emptyset$ and $\partial \delta \cap \mathbb{B}
\neq \emptyset.$ If $\partial \delta \cap A \neq \emptyset$ and
$\partial \delta \cap B \neq \emptyset$ then move the endpoint of
$\delta$ lying on $A$ to the neighboring $\alpha$-arc. Possibly
changing the orientation of $\delta$ we obtain a curve with the
required properties.

Now perform a finger move of the $\alpha$- or $\beta$-arc through
$\delta(0),$ pushing it all the way along $\delta.$ Since $R' = R
\setminus \delta$ is connected we obtain a surface diagram
$\mathcal{S}'$ where $R$ is replaced by a component homeomorphic to
$R',$ plus an extra bigon. The homeomorphism type of every other
component remains unchanged. Observe that $\chi(R') = \chi(R)+1,$ so
we have $a(\mathcal{S}') = a(\mathcal{S})-1.$ If we repeat this
process we end up in a finite number of steps with a diagram, also
denoted by $\mathcal{S},$ where $a(\mathcal{S})=0.$ Note that for
every connected surface $F$ with non-empty boundary we have $\chi(F)
\le 1,$ and $\chi(F) = 1$ if and only if $F \approx D^2.$ Thus
$a(\mathcal{S}) = 0$ implies that $R(\mathcal{S}) = \emptyset.$

Next we achieve that every component $D \in C(\mathcal{S})$ is a
bigon or a square. All the operations that follow preserve the
property that $R(\mathcal{S}) = \emptyset.$

\begin{defn} \label{defn:26}
If $D$ is a component of $\Sigma \setminus (\mathbb{A} \cup
\mathbb{B})$ then its \emph{distance} $d(D)$ from $\partial \Sigma$
is defined to be the minimum of $|\varphi \cap
(\bigcup\boldsymbol{\alpha} \cup \bigcup\boldsymbol{\beta})|$ taken
over those curves $\varphi \subset \Sigma$ for which $\varphi(0) \in
\partial \Sigma$ and $\varphi(1) \in \text{int}(D);$ furthermore, $\varphi(t)
\in \Sigma \setminus (A \cup B)$ for $0 < t \le 1.$ If $\varphi$
passes through an intersection point between an $\alpha$- and a
$\beta$-curve we count that with multiplicity two in $|\varphi \cap
(\bigcup\boldsymbol{\alpha} \cup \bigcup\boldsymbol{\beta})|.$

If $D \in C(\mathcal{S})$ is a $2n$-gon, then its \emph{badness} is
defined to be $\max \{ n-2,0 \}.$ The \emph{distance of a surface
diagram} $\mathcal{S}$ is $$d(\mathcal{S})= \max\{\,d(D) \colon D
\in C(\mathcal{S}), b(D)>0\,\}.$$

For $d > 0$ the \emph{distance $d$ complexity} of the surface
diagram $\mathcal{S}$ is defined to be the tuple
$$\left(\sum_{i=1}^m b(D_i), -b(D_1) ,\dots, -b(D_m) \right),$$
where $D_1, \dots, D_m$ are all the elements of $C(\mathcal{S})$
with $d(D) = d$ and $b(D)>0,$ enumerated such that $b(D_1) \ge \dots
\ge b(D_m).$ We order the set of distance $d$ complexities
lexicographically. Finally, let $b_d(\mathcal{S})=\sum_{i=1}^m
b(D_i).$
\end{defn}

\begin{lem} \label{lem:5}
Let $\mathcal{S}$ be a surface diagram of distance $d(\mathcal{S}) =
d > 0$ and $R(\mathcal{S}) = \emptyset.$ Then we can modify
$\mathcal{S}$ using permissible moves to get a surface diagram
$\mathcal{S}'$ with $R(\mathcal{S'}) = \emptyset,$ distance
$d(\mathcal{S}') \le d(\mathcal{S}),$ and $c_d(\mathcal{S}') <
c_d(\mathcal{S}).$
\end{lem}

\begin{proof}
Let $D_1, \dots, D_m$ be an enumeration of the distance $d$ bad
elements of $C(\mathcal{S})$ as in Definition \ref{defn:26}. Then
$D_m$ is a $2n$-gon for some $n \ge 3.$ Let $D_*$ be a component of
$\Sigma \setminus (\mathbb{A} \cup \mathbb{B})$ with $d(D_*) = d-1$
and having at least one common $\alpha$- or $\beta$-edge with $D_m.$
Without loss of generality we can suppose that they have a common
$\beta$-edge $b_*.$ Let $a_1, \dots, a_n$ be an enumeration of the
edges of $D_m$ lying in $\mathbb{A}$ starting from $b_*$ and going
around $\partial D_m$ counterclockwise.

Let $1 \le i \le n.$ We denote by $R_i^1, \dots, R_i^{k_i}$ the
following distinct components of $\Sigma \setminus (\mathbb{A} \cup
\mathbb{B}).$ For every $1 \le j \le k_i-1$ the component $R_i^j$ is
a square of distance $d(R_i^j) \ge d,$ but $R_i^{k_i}$ does not have
this property. Furthermore, $a_i \cap R_i^1 \neq \emptyset$ and
$R_i^j \cap R_i^{j+1} \subset \mathbb{A}$ for $1 \le j \le k_i-1.$
Then $R_i^{k_i}$ is either a bigon or a component of distance
$d(R_i^{k_i}) \le d.$ Note that it is possible that $R_i^{k_i} =
D_m,$ in which case $R_i^j = R_l^{k_i-j}$ for some $a_l \subset
R_i^{k_i-1} \cap R_i^{k_i}$ and every $1 \le j \le k_i-1.$

Thus if we leave $D_m$ through $a_i$ and move through opposite edges
we visit the sequence of squares $R_i^1, \dots, R_i^{k_i-1}$ until
we reach a component $R_i^{k_i}$ which is not a square of distance
$\ge d.$

Let $I = \{\,1 \le i \le n \,\colon\, R_i^{k_i} \neq D_m\,\}.$ We
claim that $I \neq \emptyset.$ Indeed, otherwise take the domain
$\mathcal{D}$ that is the sum of those components of $\Sigma
\setminus (\mathbb{A} \cup \mathbb{B})$ that appear as some $R_i^j$
for $1 \le i \le n$ and $1 \le j \le k_i,$ each taken with
coefficient one. Then $\partial \mathcal{D}$ is a sum of closed
components of $\mathbb{B}.$ Since $B$ has no closed components
$\partial \mathcal{D}$ is a sum of full $\beta$-curves,
contradicting the fact that the elements of $\boldsymbol{\beta}$ are
linearly independent in $H_1(\Sigma; \mathbb{Z}).$

First suppose that $\exists i \in I \cap \{2, \dots, n-1\}.$ Then
choose a properly embedded arc $\delta \subset D_m \cup (R_i^1 \cup
\dots \cup R_i^{k_i})$ such that $\delta(0) \in b_*$ and $\delta(1)
\in \text{int}(R_i^{k_i});$ furthermore, $|\delta \cap
\partial R_i^j| =2$ for $1 \le j < k_i.$ Observe that $\delta(t) \cap
\mathbb{B} = \emptyset$ for $0 < t \le 1.$ Do a finger move of the
$b_*$ arc along $\delta$ and call the resulting surface diagram
$\mathcal{S}'.$ The finger cuts $D_m$ into two pieces called $D_m^1$
and $D_m^2,$ and $D_*$ becomes a new component $D_*'.$

We claim that $\mathcal{S}'$ satisfies the required properties.
Indeed, $d(\mathcal{S}') \le d(\mathcal{S})$ because $\delta$ does
not enter any region of distance $<d$ except possibly $R_i^{k_i}$
for which $R_i^{k_i} \setminus \delta$ is still connected. Thus
$d(D_*') < d$ and the only new bad regions that we possibly make,
$D_m^1$ and $D_m^2,$ have a common edge with $D_*'.$ All the other
new components are bigons or squares. To show that
$c_d(\mathcal{S}') < c_d(\mathcal{S})$ we distinguish three cases.
Observe that we have
\begin{equation} \label{eqn:2}
b(D_m^1) + b(D_m^2) = b(D_m)-1.
\end{equation}
Indeed, if $D_m^1$ is a $2n_1$-gon and $D_m^2$ is a $2n_2$-gon then
$n_1 > 1$ and $n_2 > 1$ since $1 < i < n.$ Thus $b(D_m^1) = n_1 - 2$
and $b(D_m^2) = n_2 -2.$ Since the finger cuts $a_i$ into two
distinct arcs we have that $n_1 + n_2 = n + 1,$ i.e.,
$(n_1-2)+(n_2-2)=(n-2)-1.$ Furthermore, the finger cuts $R_i^j$ for
$1 \le j < k_i$ into three squares.

Case 1: $R_i^{k_i}$ is a bigon of distance $\ge d$. Then $R_i^{k_i}
\neq D_*$ because their distances are different. Thus the finger
cuts $R_i^{k_i}$ into a bigon and a square, both have badness $0.$
So equation \ref{eqn:2} implies that $b_d(\mathcal{S}') =
b_d(\mathcal{S})-1,$ showing that $c_d(\mathcal{S}') <
c_d(\mathcal{S}).$

Case 2: $d(R_i^{k_i}) < d.$ Then the finger cuts $R_i^{k_i}$ into a
bigon and a component of distance $<d.$ Thus again we have that
$b_d(\mathcal{S}') = b_d(\mathcal{S})-1.$

Case 3: $R_i^{k_i} = D_l$ for some $1 \le l < m.$ Then the finger
cuts $D_l$ into a bigon and a component $D_l'$ such that $d(D_l') =
d$ and $b(D_l') = b(D_l)+1.$ Thus $b_d(\mathcal{S}') =
b_d(\mathcal{S}).$ But we still have $c_d(\mathcal{S}') <
c_d(\mathcal{S})$ because $D_1, \dots, D_{l-1}$ remained unchanged,
$-b(D_l') < -b(D_l),$ and every other distance $d$ region in
$\mathcal{S}'$ has badness $ < b(D_l').$

Now suppose that $I \cap \{2, \dots, n-1\} = \emptyset.$ Since $I
\neq \emptyset$ we have $1 \in I$ or $n \in I.$ We can suppose
without loss of generality that $1 \in I.$ Then we have two cases.

\begin{figure}[t]
\includegraphics{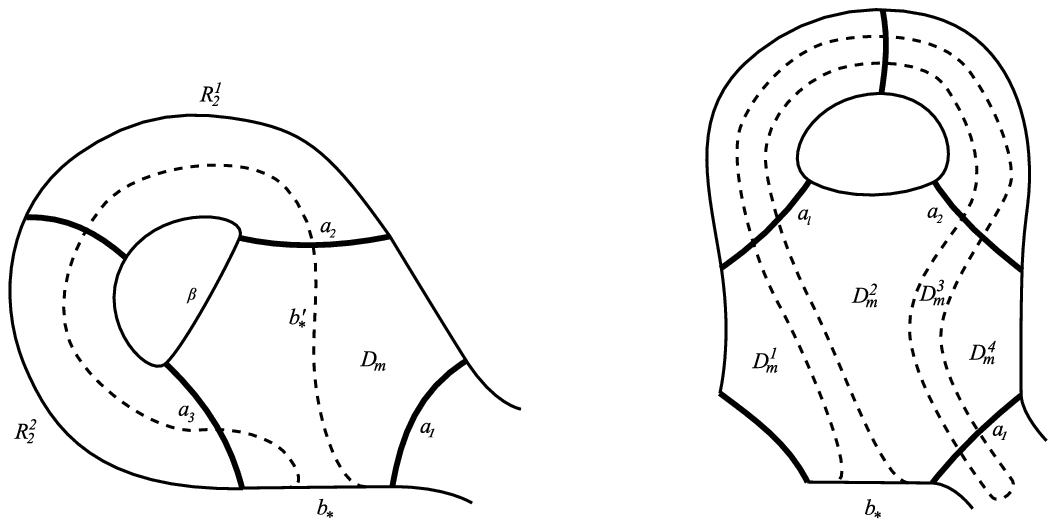}
\caption{The handle slide of Case A shown on the left. Subcase B2 is
illustrated on the right.} \label{fig:6}
\end{figure}

Case A: $n = 3;$ for an illustration see the left hand side of
Figure \ref{fig:6}. Then $R_2^{k_2} = D_m,$ and thus $R_2^{k_2-1}
\cap D_m \supset a_3,$ so $I = \{1\}.$ Let $b$ be the
$\mathbb{B}$-arc of $\partial D_m$ lying between $a_2$ and $a_3.$
Then the component $C$ of $\partial (R_2^1 \cup \dots \cup
R_2^{k_2})$ containing $b$ is a closed curve such that $C \subset
\mathbb{B}.$ Since $B$ has no closed components $C = \beta \in
\boldsymbol{\beta}$ disjoint from $b_*.$ Then handle slide $b_*$
over $\beta$ to get a new surface diagram $\mathcal{S}'.$ In
$\mathcal{S}'$ the component $D_*$ becomes $D_*'$ with $b(D_*') =
b(D_*) + 2.$ Let $b_*'$ denote $b_*$ after the handle slide. Since
$d(R_2^j) \ge d$ for $1 \le j \le k_2$ we see that $d(\mathcal{S}')
\le d(\mathcal{S});$ furthermore, $d(D_*') < d.$ The arc $b_*'$ cuts
$D_m$ into a bigon and a square; moreover, it cuts each $R_2^j$ for
$1 \le j < k_2-1$ into two squares. Thus we got rid of the distance
$d$ bad component $D_m,$ so $b_d(\mathcal{S'}) < b_d(\mathcal{S}).$

Case B: $n > 3.$ Then for some $2 < l \le n$ we have $a_l \subset
R_2^{k_2-1} \cap D_m.$

Subcase B1: $l < n;$ for an illustration see the right hand side of
Figure \ref{fig:6}. Let $$\delta \subset (R_1^1 \cup \dots \cup
R_1^{k_1}) \cup (R_2^1 \cup \dots \cup R_2^{k_2})$$ be a properly
embedded arc that starts on $b_*,$ enters $R_2^{k_2-1}$ through
$a_l,$ crosses each $R_2^j$ for $1 \le j < k_2-1$ exactly once,
reenters $D_m$ through $a_2,$ leaves $D_m$ through $a_1$ and ends in
$R_1^{k_1}.$ Note that $R_1^{k_1} \neq D_m$ since $1 \in I.$ Do a
finger move of $b_*$ along $\delta,$ we obtain a surface diagram
$\mathcal{S}'.$ The finger cuts $D_m$ into four components $D_m^1,
\dots, D_m^4$ and $D_*$ becomes a component $D_*'.$ Observe that
$D_m^3$ and $D_m^4$ are squares, $d(D_*')<d,$ and both $D_m^1$ and
$D_m^2$ have a common edge with $D_*'.$ Moreover, the only component
$\delta$ enters that can be of distance $<d$ is $R_1^{k_1}.$ Thus
$d(\mathcal{S'}) \le d(\mathcal{S}).$ Furthermore, $b(D_m^1) +
b(D_m^2) = b(D_m)-1.$ So we can conclude that $c_d(\mathcal{S'}) <
c_d(\mathcal{S})$ in a manner analogous to cases 1--3 above,
according to the type of $R_1^{k_1}$.

\begin{figure}[t]
\includegraphics{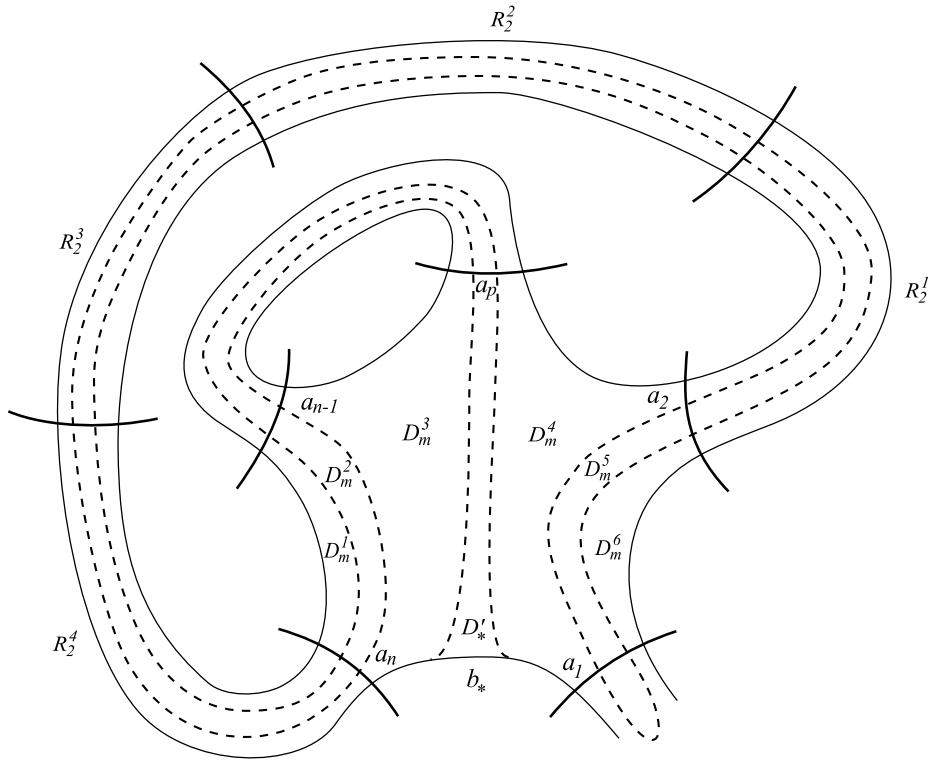}
\caption{The finger move of Subcase B2.} \label{fig:3}
\end{figure}

Subcase B2: $l =n.$ Then $a_p \subset R_{n-1}^{k_{n-1}-1} \cap D_m$
for some $2 < p< n-1.$ We define a properly embedded arc $$\delta
\subset (R_1^1 \cup \dots \cup R_1^{k_1}) \cup (R_2^1 \cup \dots
\cup R_2^{k_2}) \cup (R_p^1 \cup \dots \cup R_p^{k_p})$$ as follows
(see Figure \ref{fig:3}). The curve $\delta$ starts on $b_*,$ enters
$R_p^1$ through $a_p,$ reenters $D_m$ through $a_{n-1},$ goes into
$R_n^1 = R_2^{k_2-1}$ through $a_n,$ reenters $D_m$ through $a_2,$
leaves across $a_1,$ and ends in $R_1^{k_1}.$ Furthermore, $\delta
\cap R_i^j$ consists of a single arc for $i \in \{1,2,p\}$ and $1
\le j < k_i.$ Note that all these squares $R_i^j$ are pairwise
distinct, so $\delta$ can be chosen to be embedded. Do a finger move
of $b_*$ along $\delta$ to obtain a surface diagram $\mathcal{S}'.$
The component $D_*$ becomes $D_*'$ and the finger cuts $D_m$ into
six pieces $D_m^1, \dots, D_m^6.$ Observe that $D_m^1, D_m^2,
D_m^5,$ and $D_m^6$ are all squares; moreover, both $D_m^3$ and
$D_m^4$ have a common edge with $D_*'.$ Since $d(D_*')<d$ we have
$d(D_m^3) \le d$ and $d(D_m^4) \le d.$ Furthermore,
$b(D_m^3)+b(D_m^4) = b(D_m)-1.$ Thus we get, similarly to Subcase
B1, that $\mathcal{S}'$ has the required properties.
\end{proof}

Applying Lemma \ref{lem:5} to $\mathcal{S}$ a finite number of times
we get a surface diagram $\mathcal{S}' = (\Sigma,
\boldsymbol{\alpha}', \boldsymbol{\beta}', P)$ with $d(\mathcal{S}')
= 0,$ which means that $\mathcal{S}'$ is nice. All that remains to
show is that $(\Sigma, \boldsymbol{\alpha}', \boldsymbol{\beta}')$
is admissible if $(\Sigma, \boldsymbol{\alpha}, \boldsymbol{\beta})$
was admissible.

The proof of the fact that isotopies of the $\alpha$- and
$\beta$-curves do not spoil admissibility is a local computation
that is analogous to the one found in \cite[Section 4.3]{Sucharit}.
Handleslides only happen in Case A of Lemma \ref{lem:5}. The local
computation of \cite[Section 4.3]{Sucharit} happens in $\mathcal{D}
= R_2^1 \cup \dots \cup R_2^{k_2},$ which satisfies $\partial
\mathcal{D} \cap \mathbb{B} \subset \bigcup \boldsymbol{\beta}$
because both $b_*$ and $b$ belong to a $\beta$-curve. The
computation does not depend on whether an arc of $\partial
\mathcal{D} \cap \mathbb{A}$ belongs to $\bigcup
\boldsymbol{\alpha}$ or $B,$ so the same proof works here too.
%Alternatively, observe that if $\phi$ is a periodic domain in
%$(\Sigma, \boldsymbol{\alpha}, \boldsymbol{\beta})$ then arcs of $A$
%and $B$ appear with multiplicity zero in $\partial \phi.$ Thus if we
%consider components of $\mathbb{A}$ to be $\alpha$-curves and
%components of $\mathbb{B}$ to be $\beta$-curves then $\phi$ is a
%periodic domain in the diagram $(\Sigma, \mathbb{A}, \mathbb{B}).$
%So the local computations of \cite{Sucharit} also work in our case.

This concludes the proof of Theorem \ref{thm:7}.
\end{proof}

\section{Holomorphic disks in nice surface diagrams}

In this section we give a complete description of Maslov index one
holomorphic disks in nice balanced diagrams. Using that result we
prove Theorem \ref{thm:1}. First we state a generalization of
\cite[Corollary 4.3]{Lipshitz}.

\begin{defn} \label{defn:30}
Let $(\Sigma, \boldsymbol{\alpha}, \boldsymbol{\beta})$ be a
balanced diagram and let $\mathbf{x}, \mathbf{y} \in
\mathbb{T}_{\alpha} \cap \mathbb{T}_{\beta}.$ For $\mathcal{D} \in
D(\mathbf{x}, \mathbf{y})$ we define $\Delta(\mathcal{D})$ as
follows. Let $\phi$ be a homotopy class of Whitney disks such that
$D(\phi) = \mathcal{D}.$ Then $\Delta(\mathcal{D})$ is the algebraic
intersection number of $\phi$ and the diagonal in
$\text{Sym}^d(\Sigma).$

Suppose that $\mathcal{D} = \sum_{i=1}^m a_i\mathcal{D}_i,$ see
Definition \ref{defn:7}. If $p \in (\bigcup \boldsymbol{\alpha})
\cap (\bigcup \boldsymbol{\beta})$ and $\mathcal{D}_{i_1}, \dots,
\mathcal{D}_{i_4}$ are the four components that meet at $p$ then we
define $$n_p(\mathcal{D}) = \frac{1}{4}(a_{i_1} + \dots +
a_{i_4}).$$ Furthermore, if $\mathbf{x} = (x_1, \dots, x_d)$ and
$\mathbf{y} = (y_1, \dots, y_d)$ then let
$n_{\mathbf{x}}(\mathcal{D}) = \sum_{i=1}^d n_{x_i}(\mathcal{D})$
and $n_{\mathbf{y}}(\mathcal{D}) = \sum_{i=1}^d
n_{y_i}(\mathcal{D}).$

To define the Euler measure $e(\mathcal{D})$ of $\mathcal{D}$ choose
a metric of constant curvature $1,0,$ or $-1$ on $\Sigma$ such that
$\partial \mathcal{D}$ is geodesic and such that the corners of
$\mathcal{D}$ are right angles. Then $e(\mathcal{D})$ is $1/2\pi$
times the area of $\mathcal{D}.$
\end{defn}

\begin{rem}
The Euler measure is additive under disjoint unions and gluing of
components along boundaries. Moreover, the Euler measurer of a
$2n$-gon is $1-n/2.$
\end{rem}

\begin{prop} \label{prop:11}
If $(\Sigma, \boldsymbol{\alpha}, \boldsymbol{\beta})$ is a balanced
diagram, $\mathbf{x}, \mathbf{y} \in \mathbb{T}_{\alpha} \cap
\mathbb{T}_{\beta},$ and $\mathcal{D} \in D(\mathbf{x}, \mathbf{y})$
is a positive domain then
$$\mu(\mathcal{D}) = e(\mathcal{D}) + n_{\mathbf{x}}(\mathcal{D}) +
n_{\mathbf{y}}(\mathcal{D});$$ furthermore, $$\Delta(\mathcal{D}) =
n_{\mathbf{x}}(\mathcal{D}) + n_{\mathbf{y}}(\mathcal{D}) -
e(\mathcal{D}).$$
\end{prop}

\begin{proof}
Observe that the proof of \cite[Corollary 4.3]{Lipshitz} does not
use the fact that the number of elements of $\boldsymbol{\alpha}$
and $\boldsymbol{\beta}$ equals the genus of $\Sigma.$
\end{proof}

\begin{thm} \label{thm:8}
Suppose that $(\Sigma, \boldsymbol{\alpha}, \boldsymbol{\beta})$ is
a nice balanced diagram, $\mathbf{x}, \mathbf{y} \in
\mathbb{T}_{\alpha} \cap \mathbb{T}_{\beta},$ and $\mathcal{D} \in
D(\mathbf{x}, \mathbf{y})$ is a positive domain with
$\mu(\mathcal{D}) = 1.$ Then, for a generic almost complex
structure, $\widehat{\mathcal{M}}(\mathcal{D})$ consists of a single
element which is represented by an embedding of a disk with two or
four marked points into $\Sigma$.
\end{thm}

\begin{proof}
In light of Proposition \ref{prop:11} the proof is completely
analogous to the proofs of \cite[Theorem 3.2]{Sucharit} and
\cite[Theorem 3.3]{Sucharit}.
\end{proof}

\begin{prop} \label{prop:10}
If the surface diagram $\mathcal{S} = (\Sigma, \boldsymbol{\alpha},
\boldsymbol{\beta}, P)$ is nice and $(\Sigma, \boldsymbol{\alpha},
\boldsymbol{\beta})$ is admissible then the balanced diagram
$(\Sigma, \boldsymbol{\alpha}, \boldsymbol{\beta})$ is also nice.
\end{prop}

\begin{proof}
As before, let $C(\mathcal{S})$ denote the set of those components
of $\Sigma \setminus (\mathbb{A} \cup \mathbb{B})$ whose closure is
disjoint from $\partial \Sigma.$ Since $\mathcal{S}$ is nice each
component $R \in C(\mathcal{S})$ is a bigon or a square, and thus
its Euler measure $e(R) \ge 0.$ Let $\mathcal{S}' = (\Sigma,
\boldsymbol{\alpha}, \boldsymbol{\beta}, \emptyset).$ Then every
component $R' \in C(\mathcal{S}')$ is a sum of elements of
$C(\mathcal{S}),$ each taken with multiplicity one. Thus $e(R') \ge
0,$ which implies that $R'$ is a bigon, a square, an annulus, or a
disk. It cannot be an annulus or a disk because that would give a
nontrivial positive periodic domain in $(\Sigma,
\boldsymbol{\alpha}, \boldsymbol{\beta}).$
\end{proof}

\begin{prop} \label{prop:7}
Let $\mathcal{S} = (\Sigma, \boldsymbol{\alpha}, \boldsymbol{\beta},
P)$ be a good, nice, and admissible surface diagram and let $D(P) =
(\Sigma', \boldsymbol{\alpha}', \boldsymbol{\beta}', P_A, P_B, p).$
Then the balanced diagram $(\Sigma', \boldsymbol{\alpha}',
\boldsymbol{\beta}')$ is admissible and $$CF(\Sigma',
\boldsymbol{\alpha}', \boldsymbol{\beta}') \approx (O_P,
\partial|O_P).$$
\end{prop}

\begin{proof}
Suppose that $Q'$ is a periodic domain in $(\Sigma',
\boldsymbol{\alpha}', \boldsymbol{\beta}')$ with either no positive
or no negative multiplicities. Then $Q = p(Q')$ is a periodic domain
in $(\Sigma, \boldsymbol{\alpha}, \boldsymbol{\beta})$ since
$p(\partial Q') = \partial Q$ will be a linear combination of full
$\alpha$- and $\beta$-curves. Furthermore, $Q$ has either no
positive or no negative multiplicities, thus by the admissibility of
$(\Sigma, \boldsymbol{\alpha}, \boldsymbol{\beta})$ we get that $Q =
0.$ So $Q'$ is also zero since all of its coefficients have the same
sign.

According to Lemma \ref{lem:3} the map $p$ induces a bijection
between $\mathbb{T}_{\alpha'} \cap \mathbb{T}_{\beta'}$ and $O_P,$
which we denote by $p_*.$ We claim that $p_*$ is an isomorphism of
chain complexes.

Let $\mathbf{x}', \mathbf{y}' \in \mathbb{T}_{\alpha'} \cap
\mathbb{T}_{\beta'}$ and let $\mathbf{x} = p_*(\mathbf{x}')$ and
$\mathbf{y} = p_*(\mathbf{y}').$ Then $\mathbf{x}, \mathbf{y} \in
O_P.$ Take a positive domain $\mathcal{D}' \in D(\mathbf{x}',
\mathbf{y}')$ such that $\mu(\mathcal{D}')=1$ and let $\mathcal{D} =
p(\mathcal{D}').$ Observe that $n_{\mathbf{x}}(\mathcal{D}) =
n_{\mathbf{x}}(\mathcal{D'}),$ $n_{\mathbf{y}}(\mathcal{D}) =
n_{\mathbf{y}}(\mathcal{D'}),$ and $e(\mathcal{D}) =
e(\mathcal{D}').$ Then $\mathcal{D}$ is a positive domain with
$\mu(\mathcal{D}) = 1$ due to Proposition \ref{prop:11}. Thus $p$
induces a map $p_0$ from
$$L' = \{\,\mathcal{D}' \in D(\mathbf{x}', \mathbf{y}') \,\colon\,
\mathcal{D}' \ge 0 \,\,\text{and}\,\, \mu(\mathcal{D}') = 1 \,\}$$
to $$L = \{\,\mathcal{D} \in D(\mathbf{x}, \mathbf{y}) \,\colon\,
\mathcal{D} \ge 0 \,\,\text{and}\,\, \mu(\mathcal{D}) = 1 \,\}.$$ We
claim that $p_0$ is a bijection by constructing its inverse $r_0$.

%For $\phi' \in \mathcal{M}(\mathcal{D}')$ let $\phi = p \circ
%\phi'.$ This induces a map $p_0 \colon
%\widehat{\mathcal{M}}(\mathcal{D}') \to
%\widehat{\mathcal{M}}(\mathcal{D}).$ We show that $p_0$ is a
%bijection by constructing its inverse. According to Theorem
%\ref{thm:8} the modulus space $\widehat{\mathcal{M}}(\mathcal{D})$
%has a unique element that can be represented by an embedding $\phi =
%\phi_{\mathcal{D}} \colon D^2 \to \Sigma.$

Let $\mathcal{A} = (\cup \boldsymbol{\alpha}) \cup A$ and
$\mathcal{B} = (\cup \boldsymbol{\beta}) \cup B.$ Suppose that
$\mathcal{D} \in L.$ Then $\mathcal{D}$ is an embedded square or
bigon according to Theorem \ref{thm:8}. Let $C$ be a component of
$\mathcal{D} \cap P.$ We claim that either $\partial C \subset
\mathcal{A}$ or $\partial C \subset \mathcal{B}.$ Indeed, $C$ is a
sum of elements of $C(\mathcal{S})$ (recall that $C(\mathcal{S})$
was defined in the proof of Theorem \ref{thm:7}), which are all
bigons and squares. Thus the Euler measure $e(C) \ge 0.$ The
component $C$ cannot be an annulus or a disk since $A$ and $B$ have
no closed components and $(\Sigma, \boldsymbol{\alpha},
\boldsymbol{\beta})$ is admissible. Thus $C$ is either a bigon or a
square. Since $\mathbf{x}, \mathbf{y} \in O_P$ and because
$\mathcal{D}$ is an embedded bigon or square no corner of $C$ can be
an intersection of an $\alpha$- and a $\beta$-edge of $\partial C.$
Thus if $C$ is a bigon it can either have an $\alpha$- and an
$A$-edge, or a $\beta$- and a $B$-edge. On the other hand, if $C$ is
a square it can have two opposite $\alpha$- and two opposite
$A$-edges, or two opposite $\beta$- and two opposite $B$-edges. Note
that in all these cases if $\partial C \subset \mathcal{A}$ then
$\partial C \cap A \neq \emptyset$ and if $
\partial C \subset \mathcal{B}$ then $\partial C \cap B \neq
\emptyset.$

Now we define a map $h = h_{\mathcal{D}} \colon \mathcal{D} \to
\Sigma'$ as follows. Let $x \in \mathcal{D}.$ If $x \in \mathcal{D}
\setminus P$ then let $h(x) = p^{-1}(x).$ If $x$ lies in a component
$C$ of $\mathcal{D} \cap P$ such that $\partial C \subset
\mathcal{A}$ then let $h(x) = p^{-1}(x) \cap P_A;$ finally, let
$h(x) = p^{-1}(x) \cap P_B$ if $\partial C \subset \mathcal{B}.$ The
map $h$ is continuous because if $x \in A$ (or $x \in B$) and the
sequence $(x_n) \subset \mathcal{D} \setminus P$ converges to $x$
then the sequence $(p^{-1}(x_n))$ converges to $p^{-1}(x) \cap P_A$
(or $p^{-1}(x) \cap P_B$). See Figure \ref{fig:2}. The map $p$ is
conformal, thus $h$ is holomorphic. Furthermore, $p \circ h =
\text{id}_{\mathcal{D}}$ and thus $h$ is an embedding. So $h$ is a
conformal equivalence between $\mathcal{D}$ and $h(\mathcal{D}),$
which implies that $h(\mathcal{D}) \in L'.$ We define
$r_0(\mathcal{D})$ to be $h(\mathcal{D}).$ Then it is clear that
$p_0 \circ r_0 = \text{id}_L.$

Now we prove that $r_0 \circ p_0 = \text{id}_{L'}.$ Let
$\mathcal{D}' \in L'$ and let $\mathcal{D} = p_0(\mathcal{D}');$
furthermore, $h = h_{\mathcal{D}}.$ Since $\mathcal{D}' \ge 0$ and
$\mathcal{D}$ has only $0$ and $1$ multiplicities we see that
$\mathcal{D}'$ also has only $0$ and $1$ multiplicities. Since $p$
is conformal the map $p|\mathcal{D}' \colon \mathcal{D}' \to
\mathcal{D}$ is a conformal equivalence. Let $$h' =
(p|\mathcal{D}')^{-1} \colon \mathcal{D} \to \mathcal{D}'.$$ It
suffices to show that $h = h'$ because this would imply that
$$r_0(\mathcal{D}) = h(\mathcal{D}) = h'(\mathcal{D}) =
\mathcal{D'}.$$ Since $p \colon (\Sigma' \setminus P) \to (\Sigma
\setminus P)$ is a conformal equivalence we get that $h|(\mathcal{D}
\setminus P) = h'|(\mathcal{D} \setminus P).$ Let $C$ be a component
of $\mathcal{D} \cap P.$ Without loss of generality we can suppose
that $\partial C \subset \mathcal{A},$ and thus $\partial C \cap A
\neq \emptyset.$ Let $x \in \partial C \cap A.$ Then $h'(C)$ is
connected, so either $h'(C) \subset P_A$ or $h'(C) \subset P_B.$ But
$h'(C) \subset P_B$ cannot happen. Indeed, then we had $$h'(x) \in
p^{-1}(A) \cap P_B \subset \partial \Sigma'.$$ Moreover, the
multiplicity of $\mathcal{D}'$ at $h'(x)$ is one, but $\mathcal{D}'$
has multiplicity zero along $\partial \Sigma',$ a contradiction. So
$h'(C) \subset P_A,$ which means that $h|C = h'|C.$

Thus $p_0$ is indeed a bijection between $L'$ and $L.$ We have seen
that if $\mathcal{D}' \in L'$ and $\mathcal{D} = p_0(\mathcal{D}')$
then both $\mathcal{D}$ and $\mathcal{D}'$ are either embedded
bigons or embedded squares; moreover, $h_{\mathcal{D}}$ is a
conformal equivalence between them. In both cases
$\widehat{\mathcal{M}}(\mathcal{D})$ and
$\widehat{\mathcal{M}}(\mathcal{D}')$ have a single element.

This implies that $p_*$ is an isomorphism between the chain
complexes $(\Sigma', \boldsymbol{\alpha}', \boldsymbol{\beta}')$ and
$(O_P,\partial|O_P).$

\end{proof}

\begin{proof}[Proof of Theorem \ref{thm:1}]
According to Lemma \ref{lem:2} it is sufficient to prove the theorem
for good decomposing surfaces. Because of Proposition \ref{prop:4}
for each good decomposing surface we can find a good surface diagram
$\mathcal{S} = (\Sigma, \boldsymbol{\alpha}, \boldsymbol{\beta}, P)$
adapted to it. This surface diagram can be made admissible using
isotopies according to Proposition \ref{prop:9}. According to
Theorem \ref{thm:7} we can achieve that $\mathcal{S}$ is nice using
permissible moves, and it still defines $(M, \gamma)$ because of
Lemma \ref{lem:4}. Now Proposition \ref{prop:5} says that if $D(P) =
(\Sigma', \boldsymbol{\alpha}', \boldsymbol{\beta}', P_A, P_B, p)$
then $(\Sigma', \boldsymbol{\alpha}', \boldsymbol{\beta}')$ is a
balanced diagram defining $(M', \gamma').$ From Proposition
\ref{prop:7} we see that $(\Sigma', \boldsymbol{\alpha}',
\boldsymbol{\beta}')$ is admissible; furthermore, $$SFH(M', \gamma')
= SFH(\Sigma', \boldsymbol{\alpha}', \boldsymbol{\beta}') \approx
H(O_P,\partial|O_P).$$ Finally, Lemma \ref{lem:3} implies that
$(O_P,\partial|O_P)$ is the subcomplex of $CF(\Sigma,
\boldsymbol{\alpha}, \boldsymbol{\beta})$ generated by those
$\mathbf{x} \in \mathbb{T}_{\alpha} \cap \mathbb{T}_{\beta}$ for
which $\mathfrak{s}(\mathbf{x}) \in O_S.$ So $$H(O_P,
\partial|O_P) \approx \bigoplus_{\mathfrak{s} \in O_S} SFH(M,\gamma,
\mathfrak{s}),$$ which concludes the proof.
\end{proof}

%\section{Applications}

%First we are going to remind the reader of \cite[Definition
%2.4]{Scharlemann}, \cite[Definition 4.18]{Scharlemann}, and
%\cite[Theorem 4.19]{Scharlemann}.

%\begin{defn} \label{defn:23}
%A decomposing surface $S \subset (M, \gamma)$ is called
%\emph{conditioned} if

%(a) all arcs of $\partial S$ in any annulus component of $A(\gamma)$
%are oriented in the same direction,

%(b) no collection of simple closed curves of $\partial S \cap
%R(\gamma)$ is trivial in $H_1(R(\gamma), \partial R(\gamma)).$

%\end{defn}

%\begin{defn} \label{defn:24}
%A \emph{taut sutured manifold hierarchy} is a finite sequence
%$$(M_0, \gamma_0)\rightsquigarrow^{S_1} (M_1,\gamma_1)
%\rightsquigarrow^{S_2} \dots \rightsquigarrow^{S_n} (M_n, \gamma_n)
%$$ of taut sutured manifold decompositions for which

%(a) each $S_i$ is either a conditioned surface, a product disk, or a
%nontrivial product annulus,

%(b) no closed component of any $S_i$ separates,

%(c) $H_2(M_n, \partial M_n) = 0,$ so in fact $M_n$ is a union of
%3-balls.
%\end{defn}

%\begin{thm} \label{thm:5}
%Any taut sutured manifold $(M_0, \gamma_0)$ admits a taut sutured
%manifold hierarchy.
%\end{thm}

\section{Applications}

First we are going to remind the reader of \cite[Definition
4.1]{Gabai} and \cite[Theorem 4.2]{Gabai}. See also \cite[Theorem
4.19]{Scharlemann}.

\begin{defn} \label{defn:23}
A \emph{sutured manifold hierarchy} is a sequence of decompositions
$$(M_0, \gamma_0)\rightsquigarrow^{S_1} (M_1,\gamma_1)
\rightsquigarrow^{S_2} \dots \rightsquigarrow^{S_n} (M_n,
\gamma_n),$$ where $(M_n, \gamma_n)$ is a product sutured manifold,
i.e., $(M_n,\gamma_n) = (R \times I,
\partial R \times I)$ and $R_+(\gamma_n) = R \times \{1\}$ for some
surface $R.$ The \emph{depth} of the sutured manifold
$(M_0,\gamma_0)$ is defined to be the minimum of such $n$'s.
\end{defn}

\begin{thm} \label{thm:5}
Let $(M, \gamma)$ be a connected taut sutured manifold (see
Definition \ref{defn:37}), where $M$ is not a rational homology
sphere containing no essential tori. Then $(M, \gamma)$ has a
sutured manifold hierarchy such that each $S_i$ is connected, $S_i
\cap
\partial M_{i-1} \neq \emptyset$ if $\partial M_{i-1} \neq
\emptyset,$ and for every component $V$ of $R(\gamma_i)$ the
intersection $S_{i+1} \cap V$ is a union of parallel oriented
nonseparating simple closed curves or arcs.
\end{thm}

\begin{proof}[Proof of Theorem \ref{thm:2}]
According to Theorem \ref{thm:5} every taut balanced sutured
manifold $(M, \gamma) = (M_0, \gamma_0)$ admits a sutured manifold
hierarchy $$(M_0, \gamma_0)\rightsquigarrow^{S_1} (M_1,\gamma_1)
\rightsquigarrow^{S_2} \dots \rightsquigarrow^{S_n} (M_n,
\gamma_n).$$ Note that by definition $M$ is open. So every surface
$S_i$ in the hierarchy satisfies the requirements of Theorem
\ref{thm:1}. Thus for every $1 \le i \le n$ we get that $$SFH(M_i,
\gamma_i) \le SFH(M_{i-1}, \gamma_{i-1}).$$ Finally, since $(M_n,
\gamma_n)$ is a product it has a balanced diagram with
$\boldsymbol{\alpha} = \emptyset$ and $\boldsymbol{\beta} =
\emptyset,$ and thus $SFH(M_n, \gamma_n) \approx \mathbb{Z}$ (also
see \cite[Proposition 9.4]{sutured}). So we conclude that
$\mathbb{Z} \approx SFH(M_n,\gamma_n) \le SFH(M_0,\gamma_0).$
\end{proof}

\begin{proof}[Proof of Theorem \ref{thm:3}]
Let $Y(K)$ be the balanced sutured manifold $(M,\gamma),$ where $M$
is the knot complement $Y \setminus N(K)$ and $s(\gamma)$ consists
of a meridian of $K$ and a parallel copy of it oriented in the
opposite direction, see Definition \ref{defn:35}. Let $\xi$ be a
tangent vector field along $\partial N(K)$ pointing in the
meridional direction. Then $\xi$ lies in $v_0^{\perp},$ and thus
gives a canonical trivialization $t_0$ of $v_0^{\perp}.$ Observe
that there is a surface decomposition $$Y(K)\rightsquigarrow^S
Y(S).$$ Since $Y(S)$ is strongly balanced we can apply Theorem
\ref{thm:4} to get that
$$SFH(Y(S)) = \bigoplus_{\mathfrak{s} \in \text{Spin}^c(Y(K)) \colon
\langle\, c_1(\mathfrak{s},t_0),[S]\,\rangle = c(S,t_0)} SFH(Y(K),
\mathfrak{s}).$$ Recall that $$c(S, t_0) = \chi(S) + I(S) -
r(S,t_0).$$ Since $\partial S \subset
\partial N(K)$ is a longitude of $K$ we see that the rotation of
$p(\nu_S)$ with respect to $\xi$ is zero. Furthermore, $\chi(S) =
1-2g(S)$ and $I(S) = -1$ by part (1) of Lemma \ref{lem:9}, thus
$c(S, t_0) = -2g(S).$ So we get that
$$SFH(Y(S)) = \bigoplus_{\mathfrak{s} \in \text{Spin}^c(Y(K)) \colon
\langle\, c_1(\mathfrak{s},t_0),[S]\,\rangle = -2g(S)} SFH(Y(K),
\mathfrak{s}),$$  which in turn is isomorphic to
$\widehat{HFK}(Y,K,[S],-g(S)) \approx \widehat{HFK}(Y,K,[S],g(S)),$
see \cite{OSz2}. Note that we get $\widehat{HFK}(Y,K,[S],g(S))$ if
we decompose along $-S$ instead of $S.$
\end{proof}

Using our machinery we give a simpler proof of the fact that knot
Floer homology detects the genus of a knot, which was first proved
in \cite{OSz6}.

\begin{cor} \label{cor:2}
Let $K$ be a null-homologous knot in a rational homology 3-sphere
$Y$ whose Seifert genus is $g(K).$ Then
$$\widehat{HFK}(K,g(K)) \neq 0;$$ moreover,
$$\widehat{HFK}(K,i)=0 \,\, \text{for} \,\, i>g(K).$$
\end{cor}

\begin{proof}
First suppose that $Y \setminus N(K)$ is irreducible. Let $S$ be a
Seifert surface of $K.$ Then $Y(S)$ is taut if and only if $g(S) =
g(K).$ Thus, according to Theorem \ref{thm:2}, if $g(S) = g(K)$ then
$\mathbb{Z} \le SFH(Y(S))$ and because of \cite[Proposition
9.18]{sutured} we have that $SFH(Y(S)) = 0$ if $g(S) > g(K).$ Since
for every $i \ge g(K)$ we can find a Seifert surface $S$ such that
$g(S) = i,$ together with Theorem \ref{thm:3} we are done for the
case when $Y \setminus N(K)$ is irreducible.

Now suppose that $Y(K)$ can be written as a connected sum $(M,
\gamma) \# Y_1,$ where $(M, \gamma)$ is irreducible and $Y_1$ is a
rational homology 3-sphere. Since we can find a minimal genus
Seifert surface $S$ lying entirely in $(M, \gamma)$ (otherwise we
can do cut-and-paste along the connected sum sphere) we can apply
the connected sum formula \cite[Proposition 9.15]{sutured} to get
that $SFH(Y(S)) \approx SFH(M,\gamma) \otimes \widehat{HF}(Y_1)$
over $\mathbb{Q}.$ Since $\text{rk}\,\widehat{HF}(Y_1) \neq 0$ (see
\cite[Proposition 5.1]{OSz8}) we can finish the proof as in the
previous case.
\end{proof}

Next we are going to give a new proof of \cite[Theorem 1.1]{OSz7}.
Let $L$ be a link in $S^3,$ then $$x \colon H_2(S^3,L;\mathbb{R})
\to \mathbb{R}$$ denotes the Thurston semi-norm. Link Floer homology
provides a function $$y \colon H^1(S^3 \setminus L; \mathbb{R}) \to
\mathbb{R}$$ defined by $$y(h) = \max_{\{\,\mathfrak{s} \in
H_1(L;\mathbb{Z}) \colon \widehat{HFL}(L,\mathfrak{s}) \neq 0\,\} }
|\langle\, \mathfrak{s}, h \,\rangle|.$$

\begin{thm} \label{thm:6}
For a link $L \subset S^3$ with no trivial components and every $h
\in H^1(S^3 \setminus L)$ we have that $$2y(h) = x(PD[h]) +
\sum_{i=1}^l |\langle\, h,\mu_i \,\rangle|,$$ where $\mu_i$ is the
meridian of the $i^{\text{th}}$ component $L_i$ of $L.$
\end{thm}

\begin{proof}
Let $\xi$ be a unit vector field along $\partial N(L)$ that points
in the direction of the meridian $\mu_i$ along $\partial N(L_i).$
Consider the balanced sutured manifold $(M, \gamma) = S^3(L),$ then
$\xi$ is a section of $v_0^{\perp},$ and consequently it defines a
canonical trivialization $t_0$ of $v_0^{\perp}.$ Let $R$ be a
Thurston norm minimizing representative of $PD[h]$ having no $S^2$
components. Note that $R$ has no $D^2$ components because no
component of $L$ is trivial.

\begin{figure}[t]
\includegraphics{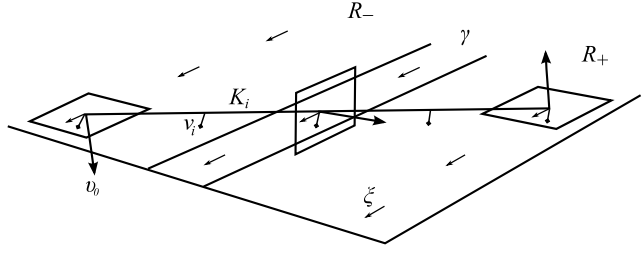}
\caption{A portion of the torus $\partial N(L_i),$ together with the
trivialization $\xi$ of $v_0^{\perp}$ and $\nu_i.$} \label{fig:9}
\end{figure}

We claim that $r(R,t_0) = 0.$ Indeed, $K_i = R \cap \partial N(L_i)$
is a torus link. We can arrange that $K_i$ and $\xi$ make a constant
angle and that $R$ is perpendicular to $\partial N(L_i)$ along
$K_i.$ Then $\nu_i = \nu_R|K_i$ is the positive unit normal field of
$K_i$ in $\partial N(L_i)$ and $\langle\, \nu_i,\xi\,\rangle_q$ is
some constant $c_i$ for every $q \in K_i,$ see Figure \ref{fig:9}.
First suppose that $c_i = 0.$ Then $K_i$ is a meridian of $L_i$ and
we can suppose that $K_i \subset R(\gamma).$ Thus $p(\nu_R)|K_i$ is
always perpendicular to $\xi.$ Now suppose that $c_i \neq 0.$  We
define the function $$a_i(q) = \langle\, p(\nu_R)/\norm{p(\nu_R)},
\xi \,\rangle_q$$ for $q \in K_i.$ Then $a_i(q) = \text{sgn}(c_i)$
for $q \in K_i \cap s(\gamma)$ and $a_i(q) = c_i$ for every $q \in
K_i \cap R(\gamma)$ such that $v_0$ is perpendicular to $R(\gamma).$
Moreover, the range of $a_i$ is $[c_i, \text{sgn}(c_i)],$ see Figure
\ref{fig:9}. So in both cases the rotation of $p(\nu_R)|K_i$ in the
trivialization $t_0$ is zero as we go around $K_i.$

%It is either a longitude of $L_i$ in which case we can arrange that
%$p(\nu_R)|K_i$ agrees with $\pm \xi$ at each point of $K_i$, or if
%the slope of $K_i$ is non-zero then $p(\nu_R)|K_i$ never agrees with
%$\xi.$ In both cases the rotation of $p(\nu_R)|K_i$ along $K_i$ with
%respect to $t_0$ is zero.

Furthermore, we can achieve that $$|\partial R \cap s(\gamma)| = 2
\sum_{i=1}^l |\langle\, h,\mu_i \,\rangle|.$$ Since $R$ is norm
minimizing and has no $S^2$ and $D^2$ components $\chi(R) =
-x(PD[h]).$ So using part (1) of Lemma \ref{lem:9} we get that
$$c(R,t_0) = -x(PD[h]) - \sum_{i=1}^l |\langle\, h,\mu_i
\,\rangle|.$$ Note that $c(R,t_0) \le 0.$

Now observe that $S^3(R)$ can be obtained from $S^3(L)$ by
decomposing along $R.$ Since $R$ is norm minimizing $S^3(R)$ is a
connected sum of taut balanced sutured manifolds, thus combining
Theorem \ref{thm:2} with the connected sum formula \cite[Proposition
9.15]{sutured} we get that $\text{rk} \, SFH(S^3(R)) \neq 0.$ So if
we apply Theorem \ref{thm:4} to the decomposition
$$S^3(L)\rightsquigarrow^R S^3(R)$$ we see that there is an
$\mathfrak{s}_0 \in\text{Spin}^c(S^3(L))$ such that $\langle\,
c_1(\mathfrak{s}_0, t_0),h \,\rangle = c(R, t_0)$ and
$\widehat{HFL}(L,\mathfrak{s}_0) \approx SFH(S^3(L),\mathfrak{s}_0)
\otimes \mathbb{Z}_2 \neq 0,$ see \cite[Proposition 9.2]{sutured}.
Thus $$2y(h) = \max_{\{\,\mathfrak{s} \in H_1(L;\mathbb{Z}) \colon
\widehat{HFL}(L,\mathfrak{s}) \neq 0\,\} } |\langle\,
c_1(\mathfrak{s},t_0), h \,\rangle| \ge x(PD[h]) + \sum_{i=1}^l
|\langle\, h,\mu_i \,\rangle|.$$

To prove that we have an equality let $\mathfrak{s}$ be a
$\text{Spin}^c$ structure on $S^3(L)$ for which $$|\langle\,
c_1(\mathfrak{s},t_0), h \,\rangle| - \left(x(PD[h]) + \sum_{i=1}^l
|\langle\, h,\mu_i \,\rangle|\right) = 2d > 0.$$ The above
difference is even because $\langle\, c_1(\mathfrak{s}_0, t_0),h
\,\rangle = c(R, t_0)$ and $$\langle c_1(\mathfrak{s},t_0) -
c_1(\mathfrak{s}_0,t_0), h \rangle = \langle
2(\mathfrak{s}-\mathfrak{s}_0), h \rangle.$$ Let $R_d$ be a Seifert
surface of $L$ obtained from $R$ by $d$ stabilizations and oriented
such that $\langle c_1(\mathfrak{s},t_0), [R_d] \rangle < 0.$
Observe that $[R_d] = \pm h,$ thus $$\langle c_1(\mathfrak{s},t_0),
[R_d] \rangle = c(R,t_0) - 2d = c(R_d,t_0),$$ which implies that
$\mathfrak{s} \in O_{R_d}.$ Now $R(S^3(R_d))$ is not Thurston norm
minimizing, thus according to \cite[Proposition 9.19]{sutured} we
have that $SFH(S^3(R_d)) = 0.$ So if we apply Theorem \ref{thm:4}
again we see that
$$\widehat{HFL}(L,\mathfrak{s}) \approx SFH(S^3(L),\mathfrak{s})
\otimes \mathbb{Z}_2 \le SFH(S^3(R_d)) \otimes \mathbb{Z}_2 =0$$ for
such an $\mathfrak{s}.$
\end{proof}

\begin{rem}
Suppose that $Y$ is an oriented 3-manifold and $L \subset Y$ is a
link such that $Y \setminus N(L)$ is irreducible. Let $x \colon
H_2(Y,L,\mathbb{R}) \to \mathbb{R}$ be the Thurston semi-norm and
for $h \in H_2(Y,L;\mathbb{R})$ let $$z(h) = \max_{\{\,\mathfrak{s}
\in \text{Spin}^c(Y,L) \colon \widehat{HFL}(Y,L,\mathfrak{s}) \neq
0\,\} } |\langle\, c_1(\mathfrak{s}), h \,\rangle|.$$ Then an
analogous proof as above gives that $$z(h) = x(h) + \sum_{i=1}^l
|\langle\, h,\mu_i \,\rangle|,$$ where $\mu_i$ is the meridian of
the $i^{\text{th}}$ component of $L.$
\end{rem}

The following proposition generalizes the horizontal decomposition
formula \cite[Theorem 3.4]{fibred}.

\begin{prop} \label{prop:8}
Let $(M,\gamma)$ be a balanced sutured manifold. Suppose that $$
(M,\gamma) \rightsquigarrow^S (M', \gamma')$$ is a decomposition
such that $S$ satisfies the requirements of Theorem \ref{thm:1},
$(M',\gamma')$ is taut, and $[S]=0$ in $H_2(M,\partial M).$ The
surface $S$ separates $(M',\gamma')$ into two parts denoted by
$(M_1, \gamma_1)$ and $(M_2,\gamma_2).$ Then $$SFH(M,\gamma) \approx
SFH(M',\gamma') \approx SFH(M_1,\gamma_1) \otimes
SFH(M_2,\gamma_2)$$ over any field $\mathbb{F}.$
\end{prop}

\begin{proof}
Since $(M',\gamma')$ is taut we can apply Theorem \ref{thm:2} to
conclude that $$SFH(M',\gamma') \neq 0.$$ Together with Theorem
\ref{thm:1} this implies that $O_S \neq \emptyset.$ Fix an element
$\mathfrak{s}_0 \in O_S.$ Then for every $\text{Spin}^c$ structure
$\mathfrak{s} \in \text{Spin}^c(M,\gamma)$ the equality $$ \langle
\, \mathfrak{s} - \mathfrak{s}_0, [S] \, \rangle = 0$$ holds since
$[S] = 0.$ Thus $\mathfrak{s} \in O_S,$ see the proof of Lemma
\ref{lem:1} and Remark \ref{rem:5}. So we get that $O_S =
\text{Spin}^c(M,\gamma),$ and thus $SFH(M',\gamma') \approx
SFH(M,\gamma).$

Now we sketch an alternative proof. Let $\mathcal{S}= (\Sigma,
\boldsymbol{\alpha},\boldsymbol{\beta},P)$ be a surface diagram
adapted to $S.$ Then $D(P) =
(\Sigma',\boldsymbol{\alpha}',\boldsymbol{\beta}',P_A,P_B,p)$ (see
Definition \ref{defn:18}) can be written as the disjoint union of
two balanced diagrams $(\Sigma_1,
\boldsymbol{\alpha}_1,\boldsymbol{\beta}_1)$ and $(\Sigma_2,
\boldsymbol{\alpha}_2,\boldsymbol{\beta}_2)$ such that $P_A \subset
\Sigma_1$ and $P_B \subset \Sigma_2.$ Let $\beta_1 \in
\boldsymbol{\beta}_1$ and $\alpha_2 \in \boldsymbol{\alpha}_2$ be
arbitrary curves. Since $\beta_1 \cap P_A = \emptyset$ and $\alpha_2
\cap P_B = \emptyset$ we get that $p(\beta_1) \cap P = \emptyset$
and $p(\alpha_2) \cap P =\emptyset.$ Furthermore, $p(\beta_1) \cap
p(\alpha_2) = \emptyset.$ Thus for the surface diagram $\mathcal{S}$
the set of inner intersection points $I_P = \emptyset.$ So Theorem
\ref{thm:1} gives that $SFH(M,\gamma) \approx SFH(M',\gamma').$

Note that $(\Sigma_i,\boldsymbol{\alpha}_i,\boldsymbol{\beta}_i)$ is
a balanced diagram of $(M_i,\gamma_i)$ for $i=1,2;$ moreover,
$$CF(\Sigma,\boldsymbol{\alpha},\boldsymbol{\beta}) \approx
CF(\Sigma_1,\boldsymbol{\alpha}_1,\boldsymbol{\beta}_1) \otimes
CF(\Sigma_2,\boldsymbol{\alpha}_2,\boldsymbol{\beta}_2).$$
\end{proof}

As a corollary of this we give a simple proof of \cite[Theorem
1.1]{Yi}. The following definition can be found in \cite{Gabai3}

\begin{defn}
The oriented surface $R \subset S^3$ is a \emph{Murasugi sum} of the
compact oriented surfaces $R_1$ and $R_2$ in $S^3$ if the following
conditions are satisfied. First, $R = R_1 \cup_E R_2,$ where $E$ is
a $2n$-gon. Furthermore, there are balls $B_1$ and $B_2$ in $S^3$
such that $R_1 \subset B_1$ and $R_2 \subset B_2,$ the intersection
$B_1 \cap B_2 = H$ is a two-sphere, $B_1 \cup B_2 = S^3,$ and $R_1
\cap H = R_2 \cap H = E.$ We also say that the knot $\partial R$ is
a Murasugi sum of the knots $\partial R_1$ and $\partial R_2.$
\end{defn}

\begin{cor} \label{cor:3}
Suppose that the knot $K \subset S^3$ is the Murasugi sum of the
knots $K_1$ and $K_2$ along some minimal genus Seifert surfaces.
Then $$\widehat{HFK}(K,g(K)) \approx \widehat{HFK}(K_1, g(K_1))
\otimes \widehat{HFK}(K_2, g(K_2))$$ over any field $\mathbb{F}.$
\end{cor}

\begin{proof}
Let $R_1$ and $R_2$ be minimal genus Seifert surfaces of $K_1$ and
$K_2,$ respectively. The Murasugi sum of $R_1$ and $R_2$ is a
minimal genus Seifert surface $R$ of $K,$ see \cite{Gabai3}. By the
definition of the Murasugi sum there is an embedded 2-sphere $H
\subset S^3$ that separates $R_1$ and $R_2$ and such that $R_1 \cap
H = R_2 \cap H$ is a $2n$-gon $E$ for some $n >0.$ Thus in the
balanced sutured manifold $S^3(R)$ the disk $D = H \setminus
\text{int}(E)$ is a separating decomposing surface that satisfies
the requirements of Theorem \ref{thm:1}. Decomposition along $D$
gives the disjoint union of $S^3(R_1)$ and $S^3(R_2),$ which is
taut. Thus, according to Proposition \ref{prop:8}, $$SFH(S^3(R))
\approx SFH(S^3(R_1)) \otimes SFH(S^3(R_2))$$ over $\mathbb{F}.$
Using Theorem \ref{thm:3} we get the required formula.
\end{proof}

\begin{lem} \label{lem:8}
Suppose that $(M,\gamma)$ is a balanced sutured manifold such that
$$H_2(M;\mathbb{Z}) = 0.$$  Let $S \subset (M,\gamma)$ be a product
annulus (see Definition \ref{defn:34}) such that at least one
component of $\partial S$ is non-zero in $H_1(R(\gamma);\mathbb{Z})$
or both components of $\partial S$ are boundary-coherent in
$R(\gamma).$ If $S$ gives a surface decomposition $(M,\gamma)
\rightsquigarrow^S (M',\gamma')$ then
$$SFH(M',\gamma') \approx SFH(M,\gamma).$$
\end{lem}

\begin{proof}
With at least one orientation of $S$ both components of $\partial S$
are boundary-coherent in $R(\gamma).$ On the other hand,
$(M',\gamma')$ does not depend on the orientation of $S.$ Thus we
can suppose that both components of $\partial S$ are
boundary-coherent.

Since $S$ is connected and $\partial S$ intersects both
$R_+(\gamma)$ and $R_-(\gamma)$ we can apply Proposition
\ref{prop:4} to get a surface diagram $(\Sigma, \boldsymbol{\alpha},
\boldsymbol{\beta}, P)$ adapted to $S.$ Here $P$ is an annulus with
one boundary component being $A$ and the other one $B.$ Thus we can
isotope all the $\alpha$- and $\beta$-curves to be disjoint from
$P,$ and so $I_P = \emptyset$ for this new diagram. The balanced
diagram $(\Sigma, \boldsymbol{\alpha}, \boldsymbol{\beta})$ is
admissible due to Proposition \ref{prop:2}. Now Lemma \ref{lem:3}
implies that for every $\mathbf{x} \in \mathbb{T}_{\alpha} \cap
\mathbb{T}_{\beta}$ we have $\mathbf{x} \in O_P$ if and only if
$\mathfrak{s}(\mathbf{x}) \in O_S.$ Consequently, $CF(\Sigma,
\boldsymbol{\alpha}, \boldsymbol{\beta}, \mathfrak{s}) = 0$ for
$\mathfrak{s} \in \text{Spin}^c(M,\gamma) \setminus O_S.$ Thus
$SFH(M,\gamma, \mathfrak{s}) = 0$ for $\mathfrak{s} \not\in O_S.$
The surface $S$ satisfies the conditions of Theorem \ref{thm:1}, and
so $SFH(M',\gamma') \approx SFH(M,\gamma).$
%Let $\mathfrak{s} \in \text{Spin}^c(M,\gamma)$ and suppose that
%$\mathfrak{s} \not\in O_S.$ For every $\mathbf{x} \in
%\mathbb{T}_{\alpha} \cap \mathbb{T}_{\beta}$ we have $\mathbf{x} \in
%O_P,$ and we have seen that this implies that $\mathfrak{s}(x) \in
%O_S.$
\end{proof}

The next proposition is an analogue of the decomposition formula for
separating product annuli proved in \cite[Theorem 4.1]{fibred} using
completely different methods.

\begin{prop} \label{prop:12}
Suppose that $(M,\gamma)$ is a balanced sutured manifold such that
$H_2(M;\mathbb{Z}) = 0.$ Let $S \subset (M,\gamma)$ be a product
annulus such that at least one component of $\partial S$ does not
bound a disk in $R(\gamma).$ Then $S$ gives a surface decomposition
$(M,\gamma) \rightsquigarrow^S (M',\gamma'),$ where $SFH(M',\gamma')
\le SFH(M,\gamma).$ If we also suppose that $S$ is separating in $M$
then $SFH(M',\gamma') \approx SFH(M,\gamma).$
\end{prop}

\begin{figure}[t]
\includegraphics{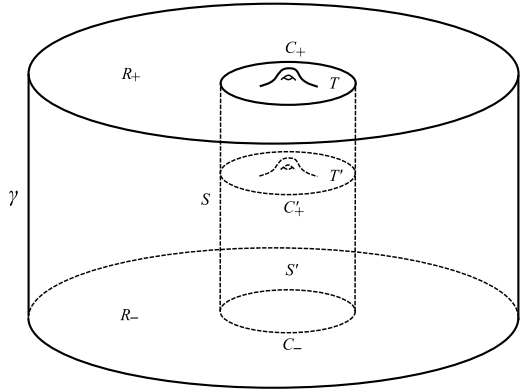}
\caption{A product annulus.} \label{fig:10}
\end{figure}

\begin{proof}
Let $C_{\pm} = \partial S \cap R_{\pm}(\gamma)$ and suppose that
$C_+$ does not bound a disk in $R_+(\gamma),$ see Figure
\ref{fig:10}. According to Lemma \ref{lem:8} we only have to
consider the case when both $[C_+]$ and $[C_-]$ are zero in
$H_1(R(\gamma);\mathbb{Z}).$ Since $(M',\gamma')$ does not depend on
the orientation of $S$ we can suppose that $S$ is oriented such that
$C_-$ is boundary-coherent in $R_-(\gamma).$ If $C_+$ is also
boundary-coherent in $R_+(\gamma)$ then we are again done due to
Lemma \ref{lem:8}. Thus suppose that $C_+$ is not boundary-coherent.

The idea of the following argument was communicated to me by Yi Ni.
Let $T$ denote the interior of $C_+$ in $R_+(\gamma);$ then $C_+$
and $\partial T$ are oriented oppositely, see Definition
\ref{defn:31}. Let $C_+'$ be a curve lying in $\text{int}(S)$
parallel and close to $C_+$ and choose a surface $T'$ parallel to
$T$ such that $\text{int}(T') \subset \text{int}(M \setminus S)$ and
$\partial T = C_+'.$ Let $S_0$ be the component of $S \setminus
C_+'$ containing $C_-.$ We define $S'$ to be the surface $S_0 \cup
T'.$ Note that the orientations of $S_0$ and $T'$ match along
$C_+',$ so $S'$ has a natural orientation. Let $(M_0,\gamma_0)$ be
the manifold obtained after decomposing $(M,\gamma)$ along $S'.$
Observe that $\partial S' = C_-$ is boundary-coherent in
$R_-(\gamma),$ thus we can apply Theorem \ref{thm:1} to $S'$ and get
that $SFH(M_0,\gamma_0) \le SFH(M,\gamma).$ If we also suppose that
$S$ is separating then $S_0$ is separating and so we have an
equality due to Proposition \ref{prop:8}. Decomposing
$(M_0,\gamma_0)$ along the annulus $S \setminus S_0$ we get a
sutured manifold homeomorphic to the disjoint union of
$(M',\gamma')$ and $(T \times I,\partial T \times I).$ Since $T \neq
D^2$ we can remove the $(T \times I,
\partial T \times I)$ part of $(M_0,\gamma_0)$ by a series of
decompositions along product disks and product annuli having no
separating boundary components. Thus $SFH(M',\gamma') \approx
SFH(M_0,\gamma_0)$ by \cite[Lemma 9.13]{sutured} and Lemma
\ref{lem:8}.
\end{proof}

\section{Fibred knots}

Ghiggini \cite{Ghiggini} (for the genus one case) and Ni
\cite{fibred} recently proved a conjecture of Ozsv\'ath and Szab\'o
that knot Floer homology detects fibred knots. We use the methods
developed in this paper to simplify their proof by avoiding the
introduction of contact structures. Moreover, we give a relationship
between knot Floer homology and the existence of depth one taut
foliations on the knot complement.

\begin{defn} \label{defn:28}
Let $(M,\gamma)$ be a balanced sutured manifold. Then $(M,\gamma)$
is called a \emph{homology product} if $H_1(M,R_+(\gamma);
\mathbb{Z}) = 0$ and $H_1(M,R_-(\gamma); \mathbb{Z}) = 0.$
Similarly, $(M,\gamma)$ is said to be a \emph{rational homology
product} if $H_1(M,R_+(\gamma); \mathbb{Q}) = 0$ and
$H_1(M,R_-(\gamma); \mathbb{Q}) = 0.$
\end{defn}

\begin{rem}
It follows from the universal coefficient theorem that every
homology product is also a rational homology product.
\end{rem}

\begin{defn} \label{defn:38}
Let $(M,\gamma)$ be a balanced sutured manifold. A decomposing
surface $S \subset M$ is called a \emph{horizontal surface} if
\begin{enumerate}[i{)}]
\item $S$ is open,
\item $\partial S \subset \gamma$ and $|\partial S|=|s(\gamma)|,$
\item $[S] = [R_+(\gamma)]$ in $H_2(M,\gamma),$
\item $\chi(S) = \chi(R_+(\gamma)).$
\end{enumerate}
We say that $(M,\gamma)$ is \emph{horizontally prime} if every
horizontal surface in $(M,\gamma)$ is parallel to either
$R_+(\gamma)$ or $R_-(\gamma).$
\end{defn}

\begin{lem} \label{lem:6}
Suppose that $(M,\gamma)$ is a balanced sutured manifold and let
$$(M,\gamma) \rightsquigarrow^S (M', \gamma')$$ be  surface
decomposition. Then the following hold.
\begin{enumerate}
\item If $(M,\gamma)$ is a rational homology product then
$H_2(M,R_{\pm}(\gamma); \mathbb{Q})=0,$ and both $H_2(M;\mathbb{Q})$
and $H_2(M; \mathbb{Z})$ vanish. \label{1}
\item If $S$ is either a product disk or a product annulus then
$(M',\gamma')$ is a rational homology product if and only if
$(M,\gamma)$ is. \label{2}
\item If $R_+(\gamma)$ is connected, $S$ is a
connected horizontal surface, and $(M,\gamma)$ is a rational
homology product then $(M',\gamma')$ is also a rational homology
product. \label{3}
\end{enumerate}
\end{lem}

\begin{proof}
Let $R_{\pm} = R_{\pm}(\gamma)$ and $R_{\pm}' = R_{\pm}(\gamma').$
Then using Alexander-Poincar\'e duality we get that $$H_2(M, R_+;
\mathbb{Q}) \approx H^1(M, R_-; \mathbb{Q}) \approx H_1(M, R_-;
\mathbb{Q}) = 0.$$ A similar argument shows that $H_2(M, R_-;
\mathbb{Q}) = 0.$

Look at the following segment of the long exact sequence of the pair
$(M,R_+):$ $$H_2(R_+; \mathbb{Q}) \to H_2(M; \mathbb{Q}) \to
H_2(M,R_+; \mathbb{Q}) = 0.$$ Since $R_+$ has no closed components
$H_2(R_+; \mathbb{Q}) = 0,$  so $H_2(M;\mathbb{Q})=0.$ From
Poincar\'e duality and the universal coefficient theorem
$$H_2(M;\mathbb{Z}) \approx H^1(M, \partial M;\mathbb{Z}) \approx
\text{Hom}(H_1(M, \partial M;\mathbb{Z}), \mathbb{Z}) \oplus
\text{Tor}(H_0(M, \partial M;\mathbb{Z})),$$ which implies that
$H_2(M;\mathbb{Z})$ is torsion free. Thus $H_2(M;\mathbb{Z}) = 0.$
This proves (\ref{1}).

Now suppose that $S$ is a product disk or a product annulus.
Consider the relative Mayer-Vietoris sequence associated to the
pairs $(M',R_+')$ and $(N(S), R_+ \cap N(S)).$ From the segment $$0
= H_1(M' \cap N(S), R_+' \cap N(S); \mathbb{Q}) \to$$ $$ \to
H_1(M',R_+'; \mathbb{Q}) \oplus H_1(N(S), R_+ \cap N(S); \mathbb{Q})
\\ \to H_1(M, R_+; \mathbb{Q}) \to$$ $$\to H_0(M' \cap N(S), R_+' \cap
N(S); \mathbb{Q})=0 $$ and since $H_1(N(S),R_+ \cap N(S);\mathbb{Q})
= 0$ we get that $H_1(M',R_+'; \mathbb{Q})=0$ if and only if
$H_1(M,R_+;\mathbb{Q}) = 0.$ We can similarly show that
$H_1(M',R_-'; \mathbb{Q}) = 0$ if and only if $H_1(M,R_-;
\mathbb{Q}) = 0.$ This proves (\ref{2}).

Finally, let $S$ be a connected horizontal surface in the balanced
sutured manifold $(M, \gamma)$ with $R_+$ connected. We denote by
$(M_1, \gamma_1)$ and $(M_2, \gamma_2)$ the two components of
$(M',\gamma'),$ indexed such that $R_+ \subset M_1$ and $R_- \subset
M_2.$ The sutured manifold $(M, \gamma)$ is a homology product and
we have already seen that this implies that $H_2(M; \mathbb{Q}) =
0.$ So from the Mayer-Vietoris sequence $$0 = H_2(S; \mathbb{Q}) \to
H_2(M_1; \mathbb{Q}) \oplus H_2(M_2; \mathbb{Q}) \to H_2(M;
\mathbb{Q}) = 0$$ we obtain that $H_2(M_i; \mathbb{Q}) = 0$ for $i
=1,2.$ Another segment of the same exact sequence is $$0 \to H_1(S;
\mathbb{Q}) \to H_1(M_1; \mathbb{Q}) \oplus H_1(M_2; \mathbb{Q}) \to
H_1(M; \mathbb{Q}) \to \widetilde{H}_0(S;\mathbb{Q}) = 0,$$ thus
$$\dim H_1(M_1;\mathbb{Q}) + \dim H_1(M_2; \mathbb{Q}) = \dim H_1(S;
\mathbb{Q}) + \dim H_1(M;\mathbb{Q}).$$ From the long exact sequence
of the pair $(M,R_{\pm})$ we see that $$0 = H_2(M,
R_{\pm};\mathbb{Q}) \to H_1(R_{\pm}; \mathbb{Q}) \to H_1(M;
\mathbb{Q}) \to 0,$$ and so $\dim H_1(M; \mathbb{Q}) = \dim
H_1(R_{\pm};\mathbb{Q}).$ Since $S$ is horizontal $\chi(S) =
\chi(R_+).$ Moreover, $R_+$ and $S$ are both connected, thus $\dim
H_1(R_{\pm};\mathbb{Q}) = \dim H_1(S;\mathbb{Q}).$ Consequently,
\begin{equation} \label{eqn:3}
\dim H_1(M_1;\mathbb{Q}) + \dim H_1(M_2;\mathbb{Q}) = 2\dim
H_1(S;\mathbb{Q}).
\end{equation}
From the long exact sequence of the triple $(M,M_2,R_-)$ consider $$
0 = H_1(M, R_-; \mathbb{Q}) \to H_1(M, M_2; \mathbb{Q}) \to
H_0(M_2,R_-;\mathbb{Q}).$$ Here $H_0(M_2,R_-;\mathbb{Q}) = 0$
because $(M_2,\gamma_2)$ is balanced. So, using excision, we get
that $H_1(M_1, S;\mathbb{Q}) \approx H_1(M, M_2;\mathbb{Q}) = 0.$
Now the exact sequence
$$0 = H_2(M_1; \mathbb{Q}) \to H_2(M_1,S; \mathbb{Q}) \to H_1(S;
\mathbb{Q}) \to H_1(M_1; \mathbb{Q}) \to H_1(M_1, S; \mathbb{Q}) =
0$$ implies that $\dim H_1(M_1;\mathbb{Q}) \le \dim H_1(S;
\mathbb{Q}).$ Using a similar argument we get that $\dim H_1(M_2;
\mathbb{Q}) \le \dim H_1(S; \mathbb{Q}).$ Together with equation
\ref{eqn:3} we see that $$\dim H_1(M_i;\mathbb{Q}) = \dim H_1(S;
\mathbb{Q})$$ for $i =1,2.$ So the map $H_1(S;\mathbb{Q}) \to
H_1(M_1;\mathbb{Q})$ is an isomorphism and we can conclude that
$H_2(M_1,S; \mathbb{Q}) = 0.$ Using Alexander-Poincar\'e duality we
get that $$H_1(M_1, R_+; \mathbb{Q}) \approx H^1(M_1, R_+;
\mathbb{Q}) \approx H_2(M_1,S; \mathbb{Q}) = 0.$$ Together with
$H_1(M_1,S;\mathbb{Q})=0$ this implies that $(M_1,\gamma_1)$ is a
rational homology product. An analogous argument shows that $(M_2,
\gamma_2)$ is also a rational homology product. This proves
(\ref{3}).
\end{proof}

Observe that the proof of \cite[Proposition 3.1]{fibred} gives the
following slightly stronger result.

\begin{lem} \label{lem:7}
Let $K$ be a null-homologous knot in the oriented 3-manifold $Y$ and
let $S$ be a Seifert surface of $K.$ If $$\text{rk}\,
\widehat{HFK}(Y,K,[S],g(S)) = 1$$ then $Y(S)$ is a homology product.
\end{lem}

\begin{cor} \label{cor:4}
If $(M,\gamma)$ is a balanced sutured manifold with $\gamma$
connected and $$\text{rk}\,SFH(M,\gamma) = 1$$ then $(M,\gamma)$ is
a homology product, and thus also a rational homology product.
\end{cor}

\begin{proof}
Since $(M,\gamma)$ is balanced and $\gamma$ is connected
$R_+(\gamma)$ and $R_-(\gamma)$ are diffeomorphic. Glue
$R_+(\gamma)$ and $R_-(\gamma)$ together using an arbitrary
diffeomorphism, then do an arbitrary Dehn filling along the torus
boundary. This way we get a 3-manifold $Y$ together with a
null-homologous knot $K$ (the core of the Dehn filling). Moreover,
$R_+(\gamma)$ gives a Seifert surface $S$ of $K$ such that $Y(S) =
(M,\gamma).$ Using Theorem \ref{thm:3}
$$\widehat{HFK}(Y,K,[S],g(S)) \approx SFH(M,\gamma).$$ So Lemma
\ref{lem:7} implies that $Y(S) = (M,\gamma)$ is a homology product.
\end{proof}

\begin{thm} \label{thm:9}
Suppose that $(M,\gamma)$ is a taut balanced sutured manifold that
is not a product. Then $SFH(M,\gamma) \ge \mathbb{Z}^2.$
\end{thm}

\begin{proof}
The outline of the proof is the following. First we modify
$(M,\gamma)$ using decompositions along product disks and product
annuli, horizontal decompositions, and adding product one-handles.
The goal is to make $(M,\gamma)$ a rational homology product,
strongly balanced, and horizontally prime. Moreover, we need a curve
in $R_+(\gamma)$ which homologically lies outside the characteristic
product region (see Definition \ref{defn:29}). Then we can find
decomposing surfaces $S_1$ and $S_2$ which give taut decompositions
$(M,\gamma) \rightsquigarrow^{S_i} (M_i,\gamma_i)$ for $i=1,2$ such
that $O_{S_1} \cap O_{S_2} = \emptyset.$ To distinguish between
$\text{Spin}^c$ structures we use Lemma \ref{lem:1}. According to
Theorem \ref{thm:2} we have $\mathbb{Z} \le SFH(M_i,\gamma_i).$ From
Theorem \ref{thm:1} we get that $$SFH(M_1,\gamma_1) \oplus
SFH(M_2,\gamma_2) \le SFH(M,\gamma),$$ which concludes the proof.

Throughout the proof we use the fact that if $(N,\nu)
\rightsquigarrow^J (N',\nu')$ is a decomposition such that $J$ is
either a product disk or product annulus then $(N,\nu)$ is taut if
and only if $(N',\nu')$ is taut. This is \cite[Lemma 3.12]{Gabai}.

By adding product one-handles to $(M,\gamma)$ as in Remark
$\ref{rem:1}$ we can achieve that $\gamma$ is connected.  This new
$(M,\gamma)$ is still taut and is not a product. It was shown in
\cite[Lemma 9.13]{sutured} that adding product one-handles does not
change $SFH(M,\gamma),$ so it is sufficient to prove the theorem
when \emph{$\gamma$ is connected}. In particular, both $R_+(\gamma)$
and $R_-(\gamma)$ are connected, thus \emph{$(M,\gamma)$ is strongly
balanced}.

%Lemma \ref{lem:7} together with Theorem \ref{thm:3} imply that if
%$(M,\gamma)$ is a balanced sutured manifold with $\gamma$ connected
%and $\text{rk}\,SFH(M,\gamma) = 1$ then $(M,\gamma)$ is a homology
%product, and thus also a rational homology product.

By Theorem \ref {thm:2} and Corollary \ref{cor:4} if the taut
balanced sutured manifold $(M,\gamma)$ is not a rational homology
product and if $\gamma$ is connected then $SFH(M,\gamma) \ge
\mathbb{Z}^2.$ So in order to prove Theorem \ref{thm:9} it is
sufficient to consider the case when \emph{$(M,\gamma)$ is a
rational homology product}.

Let $R_0, \dots, R_{k+1}$ be a maximal family of pairwise disjoint
and non-parallel horizontal surfaces in $(M,\gamma)$ such that $R_0
= R_+(\gamma)$ and $R_{k+1} = R_-(\gamma).$
%The existence of such a family follows from a theorem of Haken and Kneser.
Since $\gamma$ is connected, $R_i$ is open, and $|\partial R_i| =
|s(\gamma)|$ we get that each $R_i$ is connected. Decomposing
$(M,\gamma)$ along $R_1, \dots, R_k$ we get taut balanced sutured
manifolds $(M_i, \gamma_i)$ for $1 \le i \le k+1$ such that
$R_+(\gamma_i) = R_{i-1}$ and $R_-(\gamma_i) = R_i.$ From
Proposition \ref{prop:8}
$$SFH(M,\gamma) = \bigotimes_{i=1}^{k+1} SFH(M_i, \gamma_i)$$ over
$\mathbb{Q}.$ Furthermore, part (\ref{3}) of Lemma \ref{lem:6}
implies that each $(M_i,\gamma_i)$ is a rational homology product.
And $(M_i,\gamma_i)$ is not a product since $R_{i-1}$ and $R_i$ are
not parallel. Thus it is enough to prove Theorem \ref{thm:9} for
$(M,\gamma) =(M_1, \gamma_1).$ So we can suppose that
\emph{$(M,\gamma)$ is horizontally prime} (see Definition
\ref{defn:38}). Next we recall \cite[Definition 6.1]{fibred}, also
see \cite{cpr}.

\begin{defn} \label{defn:29}
Suppose that $(M,\gamma)$ is an irreducible sutured manifold,
$R_-(\gamma)$ and $R_+(\gamma)$ are incompressible and diffeomorphic
to each other. A \emph{product region} of $(M, \gamma)$ is a
submanifold $\Phi \times I$ of $M$ such that $\Phi$ is a compact
(possibly disconnected) surface and $\Phi \times \{0\}$ and $\Phi
\times \{1\}$ are incompressible subsurfaces of $R_-(\gamma)$ and
$R_+(\gamma),$ respectively.

In \cite[Theorem 3.4]{cpr} it is proven that there is a product
region $E \times I$ such that if $\Phi \times I$ is any product
region of $(M, \gamma)$ then there is an ambient isotopy of $M$
which takes $\Phi \times I$ into $E \times I.$ We call $E \times I$
a \emph{characteristic product region} of $(M, \gamma).$
\end{defn}

Let $E \times I$ be a characteristic product region of $(M,
\gamma).$ We can suppose that $\gamma \subset E \times I.$ Since
$(M,\gamma)$ is not a product $E \times I \neq M.$
%We claim that both $R_+(\gamma) \setminus (E \times \{1\})$ and
%$R_-(\gamma) \setminus (E \times \{0\})$ are connected. Indeed,
%otherwise there would be a horizontal surface in $(M,\gamma)$ that
%is parallel to neither $R_+(\gamma)$ nor $R_-(\gamma).$ Thus we can
%find a series of surface decompositions along product disks and
%product annuli $A$ with $\partial A$ non-separating in $R(\gamma)$
%such that in the end we get a taut balanced sutured manifold
%homeomorphic to
Let $$(M',\gamma') = (M \setminus E \times I, (\partial E \times I)
\setminus \gamma).$$ Denote the components of $(\partial E \times I)
\setminus \gamma$ by $F_1, \dots, F_m.$ Then each $F_i$ is a product
annulus in $(M,\gamma).$ Moreover, no component of $\partial F_i$
bounds a disk in $R(\gamma)$ since $E \times \{0\}$ and $E \times
\{1\}$ are incompressible subsurfaces of $R(\gamma).$ After the
sequence of decompositions along the product annuli $F_1,\dots,F_m$
we get the disjoint union of $(M',\gamma')$ and the product sutured
manifold $(E \times I, \partial E \times I).$ From part (\ref{2}) of
Lemma \ref{lem:6} we  get that $(M',\gamma')$ is also a rational
homology product. Moreover, using Proposition \ref{prop:12} and the
fact that
$$SFH((M',\gamma')\cup (E \times I, \partial E \times I)) \approx
SFH(M',\gamma') \otimes \mathbb{Z} \approx SFH(M',\gamma')$$ we
obtain that $SFH(M',\gamma') \le SFH(M,\gamma).$ Of course
$(M',\gamma')$ is not a product. Thus it is sufficient to prove that
$SFH(M', \gamma') \ge \mathbb{Z}^2.$ Note that $E' \times I =
N(\gamma')$ is a characteristic product region of $(M', \gamma').$
Furthermore, $(M',\gamma')$ is taut, horizontally prime, and
strongly balanced.

If $R_+(\gamma')$ is not planar then let $(M_1, \gamma_1) = (M',
\gamma')$ and $E_1 \times I = E' \times I.$ If $R_+(\gamma')$ is
planar then $\partial R_+(\gamma')$ is disconnected since otherwise
we had $\partial M' = S^2$ and $(M', \gamma')$ would not be
irreducible. Connect two different components of $\gamma'$ with a
product one-handle $T$ as in Remark \ref{rem:1} to obtain a sutured
manifold $(M_1,\gamma_1).$ Then $E_1 \times I = N(\gamma') \cup T$
is a characteristic product region of $(M_1, \gamma_1).$ According
to part (\ref{2}) of Lemma \ref{lem:6} the sutured manifold
$(M_1,\gamma_1)$ is also a rational homology product. In both cases
the map $$H_1(E_1 \times \{1\}; \mathbb{Q}) \to H_1(R_+(\gamma_1);
\mathbb{Q})$$ is not surjective. Indeed, in the second case the
curve $\omega$ obtained by closing the core of the handle $T \cap
R_+(\gamma_1)$ in $R_+(\gamma')$ lies outside $H_1(E_1 \times \{1\};
\mathbb{Q}).$ Also, $SFH(M_1, \gamma_1) = SFH(M', \gamma')$ in both
cases. Note that $(M_1,\gamma_1)$ is still taut, horizontally prime,
and strongly balanced.

From now on let $(M, \gamma) = (M_1, \gamma_1)$ and $E \times I =
E_1 \times I.$ Let $\omega_+ \subset R_+(\gamma)$ be a properly
embedded oriented curve such that $[\omega_+] \not\in H_1(E \times
\{1\}; \mathbb{Q}).$ Then $n[\omega_+] \not\in H_1(E \times I;
\mathbb{Z})$ for every $n \in \mathbb{Z}.$ Since $(M, \gamma)$ is a
rational homology product the maps
$$i_{\pm} \colon H_1(R_{\pm}(\gamma); \mathbb{Q}) \to H_1(M;
\mathbb{Q})$$ are isomorphism, see Lemma \ref{lem:6}. Thus there
exists a properly embedded oriented curve $\omega_- \subset
R_-(\gamma)$ such that $[\omega_-] \neq 0$ in $H_1(R_-(\gamma);
\mathbb{Q})$ and non-zero integers $a,b$ such that $a \cdot
i_+([\omega_+]) = b \cdot i_-([\omega_-])$ in $H_1(M; \mathbb{Z}).$
Choose a regular neighborhood $N(\omega_+ \cup \omega_-)$ of
$\omega_+ \cup \omega_-$ in $R(\gamma).$ Then
$$N = \gamma \cup N(\omega_+ \cup \omega_-)$$ is a subsurface of
$\partial M.$ Let $x$ be the Thurston semi-norm on
$H_2(M,N;\mathbb{Z}),$ see Definition \ref{defn:36}. Since $H_2(M;
\mathbb{Z}) = 0$ the map
$$\partial \colon H_2(M,N; \mathbb{Z}) \to H_1(N; \mathbb{Z})$$ is
injective. Thus there is a unique homology class $s \in H_2(M,N;
\mathbb{Z})$ such that $\partial s = a[\omega_+] - b[\omega_-].$
Moreover, let $$r = [R_+(\gamma)] = [R_-(\gamma)] \in H_2(M,N;
\mathbb{Z}),$$ then $\partial r = [s(\gamma)].$ We will need the
following definition, see \cite{Scharlemann}.

\begin{defn} \label{defn:32}
Suppose $(S_1,\partial S_1)$ and $(S_2, \partial S_2)$ are oriented
surfaces in general position in $(M, \partial M).$ Then the
\emph{double curve sum} of $S_1$ and $S_2$ is obtained by doing
oriented cut and paste along $S_1 \cap S_2$ to get an oriented
surface representing the cycle $S_1+S_2.$ The result in an embedded
oriented surface coinciding with $S_1 \cup S_2$ outside a regular
neighborhood of $S_1 \cap S_2.$
\end{defn}

The following claim is analogous to \cite[Lemma 6.5]{fibred}.

\begin{claim} \label{claim:1}
For any integers $p,q \ge 0$ we have a strict inequality
$$x(s+pr) + x(-s+qr) > (p+q)x(r).$$
\end{claim}

\begin{proof}
Let the surfaces $S_1$ and $S_2$ be norm minimizing representatives
of $s+pr$ and $-s+qr,$ respectively. Since $M$ is irreducible and
$R(\gamma)$ is incompressible we can assume that $S_1$ and $S_2$
have no $S^2$ or $D^2$ components. Thus $\chi(S_1) = -x(S_1)$ and
$\chi(S_2) = -x(S_2).$ Furthermore, we can suppose that $S_1$ and
$S_2$ are transversal, $(S_1\cup S_2) \cap \gamma$ consists of $p+q$
parallel copies of $s(\gamma),$ and $S_1 \cap R(\gamma) = S_2 \cap
R(\gamma)$ consists of $a$ parallel copies of $\omega_+$ and $b$
parallel copies of $\omega_-.$ Since $M$ is irreducible and $S_1$
and $S_2$ are incompressible we can achieve that $(S_1 \cup S_2)
\setminus (S_1 \cap S_2)$ has no disk components. Let $P$ denote the
double curve sum of $S_1$ and $S_2,$ see Definition \ref{defn:32}.
Then $[P] = (p+q)r$  and $P$ has no $S^2$ or $D^2$ components.
Moreover, for any double curve sum $\chi(P) = \chi(S_1)+\chi(S_2).$
Thus $x(P) = x(S_1) + x(S_2).$ Also note that $P \cap R(\gamma) =
\emptyset$ and $P \cap \gamma$ consists of $p+q$ parallel copies of
$s(\gamma).$

Suppose that $T$ is a torus component of $P.$ Then $T =
\bigcup_{j=1}^{2m} A_j,$ where $A_{2i-1} \subset S_1$ and $A_{2i}
\subset S_2$ are annuli for $1 \le i \le m.$ Let $A^1 =
\bigcup_{i=1}^m A_{2i-1}$ and $A^2 = \bigcup_{i=1}^m A_{2i},$ and
define $S_1' = (S_1 \setminus A^1) \cup (-A^2)$ and $S_2' = (S_2
\setminus A^2) \cup (-A^1).$ With a small isotopy we can achieve
that $|S_1' \cap S_2'| < |S_1 \cap S_2|.$ For $i=1,2$ we have
$\partial S_i' =
\partial S_i,$ and thus $[S_i'] = [S_i]$ in $H_2(M,N);$ moreover,
$x(S_i') = x(S_i).$ Thus we can suppose that $P$ has no torus
components.

Due to the triangle inequality we only have to exclude the case
$$x(s+pr) + x(-s+qr) = (p+q)x(r).$$ Thus suppose that $x(P) =
(p+q)x(r).$ We define a function $\varphi \colon M \setminus P \to
\mathbb{Z}$ by setting $\varphi(z)$ to be the algebraic intersection
number of $P$ with a path connecting $z$ and $R_+(\gamma).$ This is
well defined because the image of $[P] = (p+q)r$ in $H_2(M,\partial
M)$ is zero, and thus any closed curve in $M$ intersects $P$
algebraically zero times.

Let $J_i = \text{cl}(\varphi^{-1}(i))$ for $0 \le i \le p+q$ and let
$P_i = J_{i-1} \cap J_i$ for $1 \le i \le p+q.$ Then $P =
\coprod_{i=1}^{p+q} P_i$ and $\bigcup_{k=0}^{i-1} J_i$ is a homology
between $R_+(\gamma)$ and $P_i$ in $H_2(M,N).$ Thus $[P_i] =
[R_+(\gamma)] = r$ and $x(P_i) \ge x(r).$ Since $$\sum_{i=1}^{p+q}
x(P_i) = x(P) = (p+q)x(r)$$ we must have $x(P_i) = x(r)$ for $1 \le
i \le p+q.$ Each $P_i$ is connected since it has no $S^2$ and $T^2$
components, and $H_2(M)=0$ implies that $P_i$ can have no higher
genus closed components, otherwise it would not be norm minimizing
in $r.$

So each $P_i$ is a horizontal surface in $(M,\gamma),$ consequently
it is parallel to $R_+(\gamma)$ or $R_-(\gamma).$ Thus for some $0
\le k \le p+q$ the surfaces $P_1, \dots, P_k$ are parallel to
$R_+(\gamma)$ and $P_{k+1}, \dots, P_{p+q}$ are parallel to
$R_-(\gamma).$ Let $P_0 = R_+(\gamma)$ and $P_{p+q+1} =
R_-(\gamma).$

We can isotope $S_1$ such that $S_1 \cap \text{int}(J_i)$ is a
collection of vertical annuli for $0 \le i \le p+q.$ Thus $S_1 \cap
\text{int}(J_i) = C_i \times (0,1),$ where $C_i$ is a collection of
circles in $P_i.$ Let $\gamma_k = \gamma \cap J_k.$ Observe that
there is a homeomorphism $h \colon (M, \gamma) \to (J_k, \gamma_k)$
such that $[C_k] = a[h(\omega_+)]$ in $H_1(P_k).$ Since
$a[h(\omega_+)] \not \in H_1(h(E \times \{1\}))$ there is a
component $C_k'$ of $C_k$ such that $[C_k'] \not \in H_1(h(E \times
\{1\})).$ Thus the product annulus $C_k' \times I$ cannot be
homotoped into $h(E \times I),$ which contradicts the fact that $h(E
\times I)$ is a characteristic product region of $(J_k, \gamma_k).$
\end{proof}

From \cite[Theorem 2.5]{Scharlemann} we see that there are
decomposing surfaces $S_1$ and $S_2$ in $(M, \gamma)$ such that
\begin{enumerate}
\item $[S_1] = s + pr$ and $[S_2] = -s + qr$ in $H_2(M,N)$ for some
integers $p,q \ge 0,$
\item if we decompose $(M,\gamma)$ along $S_i$
for $i =1,2$ we get a \emph{taut} sutured manifold $(M_i,
\gamma_i),$
\item $\nu_{S_i}$ is nowhere parallel to $v_0$ along $\partial S_i$
for $i=1,2,$
\item $\partial S_1 \cap R(\gamma)$ consists of $a$ parallel copies of
$\omega_+$ and $b$ parallel copies of $-\omega_-,$
\item $\partial S_2 \cap R(\gamma) = -\partial S_1 \cap R(\gamma),$
\item $\partial S_i \cap \gamma$ consists of parallel copies of
$s(\gamma)$ and $\nu_{S_i}|(\partial S_i \cap \gamma)$ points out of
$M$ for $i=1,2.$
\end{enumerate}
From (2) and Theorem \ref{thm:2} we get that
$$\mathbb{Z} \le SFH(M_i,\gamma_i)$$ for $i =1,2.$ Since $(M,
\gamma)$ is strongly balanced and $S$ satisfies (3) we can define
$c(S_1,t)$ and $c(S_2,t)$ for some trivialization $t$ of
$v_0^{\perp},$ see Definition \ref{defn:14}.

Using part (2) of Lemma \ref{lem:9} and (6) we get that $I(S_1)=0$
and $I(S_2) = 0.$ Moreover, $r(S_1,t) = p\chi(R_+(\gamma)) + K$ and
$r(S_2,t) = q\chi(R_+(\gamma)) - K,$ where $K$ is the contribution
of $\partial S_1 \cap R(\gamma)$ to $r(S_1,t).$

Since $(M,\gamma)$ is taut $\chi(R_+(\gamma)) = -x(r).$ Thus
$$c(S_1,t) = \chi(S_1) + px(r) - K = -x(s+pr) + px(r) - K$$ and
$$c(S_2,t) = \chi(S_2) + qx(r) + K = -x(-s + qr) + qx(r) + K.$$ From
Claim \ref{claim:1} we get that
$$c(S_1,t) + c(S_2,t) = (p+q)x(r) - (x(s+pr) + x(-s + qr)) < 0.$$
Let $\mathfrak{s}_i \in O_{S_i}$ for $i = 1,2.$ Lemma \ref{lem:1}
implies that $\langle c_1(\mathfrak{s}_1,t),[S_1] \rangle =
c(S_1,t)$ and $\langle c_1(\mathfrak{s}_2,t), [S_2] \rangle =
c(S_2,t).$ But $r=0$ in $H_2(M,\partial M),$ and thus $[S_1]=s =
-[S_2]$ in $H_2(M,
\partial M).$ So $\langle c_1(\mathfrak{s}_2, t), [S_1] \rangle =
-c(S_2,t).$ Together with $c(S_1,t) \neq -c(S_2,t)$ this implies
that $\mathfrak{s}_1 \neq \mathfrak{s}_2,$ and thus $O_{S_1} \cap
O_{S_2} = \emptyset.$ Using Theorem \ref{thm:1} we get that
$$\mathbb{Z}^2 \le SFH(M_1,\gamma_1) \oplus SFH(M_2,\gamma_2) \le
SFH(M,\gamma).$$ This concludes the proof of Theorem \ref{thm:9}.
\end{proof}

\begin{thm} \label{thm:10}
Let $K$ be a null-homologous knot in an oriented 3-manifold $Y$ such
that $Y \setminus K$ is irreducible and let $S$ be a Seifert surface
of $K.$ If $$\text{rk} \, \widehat{HFK}(Y,K,[S],g(S)) = 1$$ then $K$
is fibred with fibre $S$.
\end{thm}

\begin{proof}
From Theorem \ref{thm:3}
$$SFH(Y(S)) \approx \widehat{HFK}(Y,K,[S],g(S)).$$ Consequently,
$SFH(Y(S)) \neq 0$ and thus $Y(S)$ is taut. So we can apply Theorem
\ref{thm:9} to $Y(S)$ and conclude that $Y(S)$ is a product, since
otherwise we had $\mathbb{Z}^2 \le SFH(Y(S)).$ This implies that the
knot $K$ is fibred with fibre $S$.
\end{proof}

\begin{thm} \label{thm:11}
Let $(M,\gamma)$ be a taut balanced sutured manifold that is a
rational homology product. If $\text{rk} \, SFH(M,\gamma) < 4$ then
the depth of $(M,\gamma)$ is at most one.
\end{thm}

\begin{proof}
Suppose that the depth of $(M,\gamma)$ is $\ge 2.$ Note that
decompositions along product disks and product annuli do not
decrease the depth of a sutured manifold. Thus applying the same
procedure to $(M,\gamma)$ as in the proof of Theorem \ref{thm:9} we
get two depth $\ge 1$ (i.e., non-product) taut balanced sutured
manifolds $(M_1, \gamma_1)$ and $(M_2,\gamma_2)$ such that
$$SFH(M,\gamma) \ge SFH(M_1,\gamma_1) \oplus SFH(M_2,\gamma_2).$$
From Theorem \ref{thm:9} we see that $SFH(M_i,\gamma_i) \ge
\mathbb{Z}^2$ for $i = 1,2.$ Thus $SFH(M,\gamma) \ge \mathbb{Z}^4.$
\end{proof}

\begin{proof}[Proof of Theorem \ref{thm:12}]
Let $S$ be a genus $g$ Seifert surface of $K.$ Then $(M,\gamma) =
Y(S)$ is a taut balanced sutured manifold with $SFH(Y(S)) \approx
\widehat{HFK}(Y,K,g)$ due to Theorem \ref{thm:3}. The linking matrix
$V$ of $S$ is a matrix of the map $$i_+ \colon
H_1(R_+(\gamma);\mathbb{Q}) \to H_1(M; \mathbb{Q}),$$ thus $\det V =
\pm a_g \neq 0$ and $i_+$ is an isomorphism. From the long exact
sequence of the pair $(M,R_+(\gamma))$ we see that
$H_1(M,R_+(\gamma); \mathbb{Q}) = 0.$ Similarly, $H_1(M,R_-(\gamma);
\mathbb{Q})$ is also zero, thus $(M,\gamma)$ is a rational homology
product. Using Theorem \ref{thm:11} we conclude that the depth of
$(M,\gamma)$ is $\le 1.$ Now using \cite{Gabai} we get a depth $\le
1$ taut foliation on $(M,\gamma)$ transverse to $\gamma$ and leaves
including $R_{\pm}(\gamma).$
\end{proof}

\begin{rem} \label{rem:6}
If $\text{rk} \, \widehat{HFK}(Y,K,g) = 3$ then using the fact that
$\chi\left(\widehat{HFK}(Y,K,g)\right) = a_g$ we see that the
condition $a_g \neq 0$ is automatically satisfied.
\end{rem}

\begin{qn}
Let $K$ be a knot in a rational homology 3-sphere $Y$ and suppose
that $k$ is a positive integer. Does
$$\text{rk} \, \widehat{HFK}(Y,K,g(K)) < 2^k$$ imply that $Y
\setminus N(K)$ has a depth $< k$ taut foliation transverse to
$\partial N(K) ?$
\end{qn}

% ----------------------------------------------------------------
\bibliographystyle{amsplain}
\bibliography{topology}
\end{document}